\documentclass[preprint,3p,sort&compress,final,times]{elsarticle}
\usepackage{amsmath}
\usepackage{amssymb}
\usepackage{mathtools}
\usepackage{placeins}
\usepackage[usenames,dvipsnames]{color}
\usepackage{xcolor}
\usepackage{cases}
\usepackage{array}
\usepackage{todonotes}
\usepackage{lineno}

\usepackage[usenames,dvipsnames]{color}

\newcommand*{\ldblbrace}{\{\mskip-5mu\{}
\newcommand*{\rdblbrace}{\}\mskip-5mu\}}

\journal{}
\begin{document}
\begin{frontmatter}

\title{Helmholtz preconditioning for the compressible Euler equations using mixed finite elements with Lorenz staggering}
\author[BOM]{David Lee\corref{cor}}
\ead{davelee2804@gmail.com}
\author[ANU]{Alberto F. Mart\'in}
\author[ANU2]{Kieran Ricardo}

\cortext[cor]{Corresponding author.}
\address[BOM]{Bureau of Meteorology, Melbourne, Australia}
\address[ANU]{School of Computing, Australian National University, Canberra, Australia.}
\address[ANU2]{Mathematical Sciences Institute, Australian National University, Canberra, Australia}

\begin{abstract}
Implicit solvers for atmospheric models are often accelerated via the solution of a preconditioned system. For block 
preconditioners this typically involves the factorisation of the (approximate) Jacobian resulting from linearization 
of the coupled system into a Helmholtz equation for some function of the pressure. 
Here we present a preconditioner for the compressible Euler 
equations with a flux form representation of the potential temperature on the Lorenz grid using mixed finite elements.
This formulation allows for spatial discretisations that conserve both energy and potential temperature variance.
By introducing the dry thermodynamic entropy as an auxiliary variable for the solution of the algebraic system, 
the resulting preconditioner is shown to have a similar block structure to an existing preconditioner for
the material form transport of potential temperature on the Charney-Phillips grid. This new formulation is also shown to be 
more efficient and stable than both the material form transport of potential temperature on the Charney-Phillips grid, and a previous
Helmholtz preconditioner for the flux form transport of density weighted potential temperature on the Lorenz grid
for a 1D thermal bubble configuration.
The new preconditioner is further verified against standard two dimensional test cases
in a vertical slice geometry.
\end{abstract}

\end{frontmatter}

\section{Introduction}

Models of compressible atmospheric dynamics must be stable for processes across a range of temporal scales, including
fast acoustic and inertio-gravity waves as well as sub-critical inertial motions. Implicit time stepping methods are often used
in order to ensure numerical stability over long time steps 
that do not explicitly resolve fast processes. The solution of the
resulting algebraic system is often preconditioned via the approximate block factorisation of the Jacobian operator for the coupled system into 
a single Helmholtz equation for some form of the pressure. This may be applied as part of a fully implicit three
dimensional solver \cite{Yeh(2002)},\cite{Wood(2014)},\cite{Melvin(2019)},\cite{Maynard(2020)}, or as part of a 
horizontally explicit-vertically implicit dimensional splitting \cite{Lee(2021)},\cite{LP(2021)},\cite{Reddy(2023)}.
If the velocity space is hybridized such that the continuity of winds across cell boundaries is enforced via 
an additional set of Lagrange multipliers rather than the unique solution of the wind at the boundaries, then 
the Helmholtz equation may alternatively be expressed in terms of these Lagrange multipliers instead of the pressure 
variable \cite{Gibson(2020)},\cite{Betteridge(2023)}.

The Jacobian matrix from which the Helmholtz operator is derived is typically approximated so as to account for the 
stiff terms associated with fast time scales, while omitting or linearising those associated with other processes. 
This is due to the computational expense of assembling the full Jacobian, and also in the case of mixed finite element
methods because the resulting Helmholtz operator is only semi-definite due to the approximate lumping of mass matrix 
inverses during block factorisation and so cannot be reliably solved to convergence.

There are numerous choices as to how the approximations to the Jacobian and the resulting Helmholtz operator are 
formulated, in addition to the choices one may make in terms of equation form and spatial discretisation. This 
article compares several different formulations of the Helmholtz operator for a dry compressible atmosphere, 
that result from different vertical placements of the potential temperature (the Lorenz and Charney-Phillips 
grids \cite{Arakawa(1996)},\cite{TW(2005)}) and different forms of the potential temperature transport equation 
(flux form and material form) within mixed finite element spatial discretisations 
\cite{Melvin(2019)},\cite{Bendall(2020)},\cite{Lee(2021)}.

There are different benefits and drawbacks to these different modelling choices. With regards
to the choice of function space for the thermodynamic variable, collocating this with the vertical velocity 
ensures that there are no spurious computational modes, as is the case when this is collocated with the pressure \cite{Arakawa(1996)},
and also provides for an optimal representation of buoyancy modes \cite{TW(2005)},\cite{Melvin(2018)}. However,
whilst collocating the thermodynamic variable with the pressure variable on 
finite element spaces with discontinuities across element boundaries
does admit a spurious buoyancy mode, it also allows for energy conserving discretisations
\cite{Taylor(2020)},\cite{Lee(2021)},\cite{LP(2021)}. This is because the variational derivative of the thermodynamic 
variable is a function of the pressure, and so collocating these ensures that the chain rule required for the 
balance of kinetic and internal energy exchanges may be discretely preserved. Note that energy conservation does
not require the collocation of all prognostic variables, only the collocation of each prognostic variable with the variational
derivative of the energy with respect to that variable. There are numerous examples of mixed finite element formulations which 
preserve this property, while still staggering the different prognostic variables, for example 
\cite{Bauer(2018)},\cite{Eldred(2019)},\cite{LP(2021)}.

If the density and the thermodynamic variable are staggered, then the transport equations for these variables
are not consistent if the thermodynamic variable is represented as a density weighted quantity in flux form
\cite{Eldred(2019)}. 
However consistent material transport may potentially still be recovered 
for such a staggering via a judicious choice of numerical fluxes, as has previously been described for the
conservative material transport of the thermodynamic variable \cite{Thuburn(2022)}, and the moisture mixing ratios \cite{Bendall(2023)}.

If the density and thermodynamic variables are collocated however, then a material form transport equation for the thermodynamic
variable may be consistently derived from its density weighted flux form equivalent in a straightforward manner for an energy conserving 
spatial discretisation. 
Consequently we are free 
to choose between flux form and material form representations of the thermodynamic variable if this is 
collocated with the density. One benefit of using flux form transport of the thermodynamic variable 
is that in this form the quadratic variance of the thermodynamic tracer constitutes a conserved mathematical entropy
of the dry compressible Euler equations. For a careful choice of numerical fluxes this invariant 
can be discretely conserved, which can improve model stability for grid scale thermal oscillations \cite{Ricardo(2024)}.

In the present article we compare four choices of mixed finite element formulations of the compressible Euler 
equations from those referenced above, in terms of the stability and convergence of the Helmholtz operators implied 
by their discretisation. These include Helmholtz operators for the material transport of potential temperature, 
collocated with the vertical velocity (the \emph{Charney-Phillips} grid) \cite{Melvin(2019)}; the flux form 
transport of density weighted potential temperature collocated with the Exner pressure (the \emph{Lorenz} grid) 
\cite{Lee(2021)}; the material transport of potential temperature on the Lorenz grid; and a novel alternative formulation of the flux form potential temperature on the Lorenz grid 
where this is re-scaled to derive a material form transport equation for the thermodynamic entropy, leading to a block structure 
and temporal scalings that more closely match those for material transport on the Charney-Phillips grid.

The remainder of this article is structured as follows. In Section~\ref{sec::mixed_fem} the mixed finite element discretisation of the 
compressible Euler equations will be introduced. For a more detailed discussion of this subject the reader is 
directed to the references therein. Section~\ref{sec::helmholtz} details the formulation of a novel Helmholtz preconditioner for 
the flux form representation of the density weighted potential temperature on the Lorenz grid, and compares this 
to previous Helmholtz operator formulations. Results comparing these four different formulations are presented in Section 
~\ref{sec::results}, including both detailed comparisons in 1D, and reproductions of standard test cases for the new 
preconditioner in 2D. Finally, conclusions and future research directions will be discussed in Section~\ref{sec::conclusions}.

\section{Mixed finite element discretisation of the compressible Euler equations}
\label{sec::mixed_fem}

The compressible Euler equations for atmospheric motion may be described in terms of the
velocity, $\boldsymbol{u}$, density, $\rho$, density weighted potential temperature, $\Theta$,
and the Exner pressure, $\Pi$ as
\begin{subequations}\label{eq::ce_cont}
\begin{align}
	\frac{\partial\boldsymbol{u}}{\partial t} + \boldsymbol{q}\times\boldsymbol{F} + 
	\nabla\Phi + \theta\nabla\Pi &= 0\\
	\frac{\partial\rho}{\partial t} + \nabla\cdot\boldsymbol{F} &= 0\label{eq::rho_cont}\\
	\frac{\partial\Theta}{\partial t} + \nabla\cdot(\theta\boldsymbol{F}) &= 0\label{eq::Theta_cont}\\
	\Pi - c_p\Big(\frac{R\Theta}{p_0}\Big)^{R/c_v} &= 0 \ ,\label{eq::Pi_cont}
\end{align}
\end{subequations}
where $p_0$ is the reference pressure, $R=c_p-c_v$ is the ideal gas constant, $\boldsymbol{F}=\rho\boldsymbol{u}$ 
is the momentum, $\Phi=\frac{1}{2}\boldsymbol{u}\cdot\boldsymbol{u} + gz$ is the Bernoulli potential (with $z$
being the vertical coordinate), 
$\theta=\Theta/\rho$ is the potential temperature and
\begin{equation}
	\boldsymbol{q}=\frac{\nabla\times\boldsymbol{u} + f\hat{\boldsymbol{k}}}{\rho}\ ,
\end{equation}
is the potential vorticity, $f$ is the Coriolis term, $g$ is the gravitational constant 
and $\hat{\boldsymbol{k}}$ is the unit vector radially pointing outwards.

This system of partial differential equations may be spatially discretised using a mixed finite element method 
\cite{Natale(2016)},\cite{Melvin(2019)},\cite{Lee(2021)}. To this end, we first introduce the discrete 
subspaces over the three dimensional domain $\Omega$ as $\mathbb{W}_0\subset H_1(\Omega)$, 
$\mathbb{W}_1\subset H(\mathrm{curl},\Omega)$, $\mathbb{W}_2\subset H(\mathrm{div},\Omega)$,
$\mathbb{W}_3\subset L^2(\Omega)$. The Lorenz staggering of the solution variables may be 
achieved for the mixed finite element form of the discrete problem by assigning the
discrete solution variables (denoted by the subscript $h$) to the finite element spaces as 
$\boldsymbol{u}_h\in\mathbb{W}_2$, $\rho_h,\Theta_h,\Pi_h\in\mathbb{W}_3$. Additionally, we 
may assign the diagnostic fields as $\boldsymbol{q}_h\in\mathbb{W}_1$, $\boldsymbol{F}_h\in\mathbb{W}_2$,
$\Phi_h,\theta_h\in\mathbb{W}_3$. 

We elect to use the flux form of the potential temperature, 
\eqref{eq::Theta_cont}, and represent this discretely on the Lorenz grid ($\Theta_h\in\mathbb{W}_3$),
since doing so ensures the conservation of energy, as discussed below, as well as the consistent 
material transport of $\theta_h\in\mathbb{W}_3$ \cite{Eldred(2019)}. 
As discussed in the introduction, this choice can allow for discretisations that conserve 
potential temperature variance, and therefore improve model stability for fine scale thermal oscillations
\cite{Ricardo(2024)}.

The discrete variational form of the system is then derived by multiplying the equations in \eqref{eq::ce_cont} by the test functions $\boldsymbol{v}_h\in\mathbb{W}_2$,
$\phi_h,\psi_h,\chi_h\in\mathbb{W}_3$, respectively, and integrating over the domain $\Omega$ for a time step
of $\Delta t$ between time level $n$ and the new time level at nonlinear iteration $k$ as

\begin{subequations}\label{eq::ce_disc}
\begin{align}
	\int\boldsymbol{v}_h\cdot(\boldsymbol{u_h}^k-\boldsymbol{u}_h^n) + \label{eq::ce_disc_u}
	\Delta t\boldsymbol{v}_h\cdot\overline{\boldsymbol{q}}_h\times\overline{\boldsymbol{F}}_h - 
	\Delta t(\nabla\cdot\boldsymbol{v}_h)\overline{\Phi}_h - 
	\Delta t\nabla\cdot(\boldsymbol{v}_h\overline{\theta}_h)\overline{\Pi}_h\mathrm{d}\Omega& \\\notag
	+\Delta t\int[\![\overline{\theta}_h\boldsymbol{v}_h\cdot\hat{\boldsymbol{n}}_{\Gamma}]\!]\ldblbrace\overline{\Pi}_h\rdblbrace - 
	c[\![\overline{\theta}_h\boldsymbol{v}_h\cdot\hat{\boldsymbol{n}}_{\Gamma}]\!][\![\overline{\Pi}_h]\!]\mathrm{d}\Gamma
	&= \mathcal{R}_u\\
	\int\phi_h(\rho_h^k-\rho_h^n) + \Delta t\phi_h\nabla\cdot\overline{\boldsymbol{F}}_h\mathrm{d}\Omega& = \mathcal{R}_{\rho} \label{eq::ce_disc_rho}\\
	\int\psi_h(\Theta_h^k-\Theta_h^n) + \label{eq::ce_disc_Theta}
	\Delta t\psi_h\nabla\cdot(\overline{\theta}_h\overline{\boldsymbol{F}}_h)\mathrm{d}\Omega -
	\Delta t\int[\![\overline{\theta}_h\overline{\boldsymbol{F}}_h\cdot\hat{\boldsymbol{n}}_{\Gamma}]\!]\ldblbrace\psi_h\rdblbrace -
	c[\![\overline{\theta}_h\overline{\boldsymbol{F}}_h\cdot\hat{\boldsymbol{n}}_{\Gamma}]\!][\![\psi_h]\!]\mathrm{d}\Gamma 
	&= \mathcal{R}_{\Theta}\\
	\int\chi_h\log(\Pi_h^k) - \chi_h\frac{R}{c_v}\log(\Theta_h^k) - \label{eq::ce_disc_Pi}
	\chi_h\frac{R}{c_v}\log\Big(\frac{R}{p_0}\Big) - \chi_h\log(c_p)\mathrm{d}\Omega &= \mathcal{R}_{\Pi}\ ,
\end{align}
\end{subequations}
\noindent
where $\mathcal{R}_u$, $\mathcal{R}_{\rho}$, $\mathcal{R}_{\Theta}$, $\mathcal{R}_{\Pi}$ are the different 
components of the discrete nonlinear residual functional, and
$c = sign(\overline{\boldsymbol{F}}_h\cdot\hat{\boldsymbol{n}}_{\Gamma})/2$ in \eqref{eq::ce_disc_Theta} 
indicates an upwinded density weighted potential temperature flux, as well as its energy conserving adjoint in 
\eqref{eq::ce_disc_u}, and $c=0$ indicates a centered flux.
Note that in \eqref{eq::ce_disc_Pi} we have taken the natural logarithm of the original equation of state 
\eqref{eq::Pi_cont}. As in \cite{Lee(2021)},\cite{LP(2021)} this is done in order to represent the nonlinearity of the 
equation of state more smoothly than as a fractional exponent in its original form, and also to ensure a simpler 
linearisation of these terms within the approximate Jacobian operator as given in Section 3.
Note also that integration by parts has been selectively applied where the trial space does not support
the necessary differential operator. We have also included internal boundary integrals over the 
element faces, $\Gamma$, for which the unit normal vector is defined as $\hat{\boldsymbol{n}}_{\Gamma}$, 
involving both jumps, $[\![a_h]\!]:=a_h^+-a_h^-$ and means, $\ldblbrace a\rdblbrace:=(a_h^++a_h^-)/2$, where 
$a_h^+$ and $a_h^-$ refer to the evaluation of a discrete quantity $a_h$ from the elements on the right and 
left side of the face respectively. 

The energy conserving Exner pressure gradient and potential temperature transport terms in \eqref{eq::ce_disc_u}
and \eqref{eq::ce_disc_Theta} respectively are conceptually similar to those previously described for the thermal
shallow water equations \cite{Eldred(2019)}, however in the present case we apply a strong form of the potential temperature
transport term and a weak form of the Exner pressure gradient term. In contrast, \cite{Eldred(2019)} uses a weak form transport 
term and strong form gradient term, implying that the signs of the boundary fluxes are opposite in these two different formulations.

The overbars $\overline{a}_h$ in \eqref{eq::ce_disc} denotes exact second order time integration of $a_h$. For linear quantities
this is nothing more than a simple averaging. The exception to this are the quadratic nonlinear terms
$\overline{\boldsymbol{F}}_h$ and $\overline{\Phi}_h$, for which this is computed using Simpson's rule. 
Doing so ensures that these are exactly integrated in time for a second order temporal integration,
which is necessary for exact temporal conservation of energy, as previously shown for the shallow water
\cite{Bauer(2018)}, thermal shallow water \cite{Eldred(2019)} and compressible Euler equations \cite{Lee(2021)}. 
For the compressible Euler equations the total energy is given as
\begin{equation}\label{eq::energy}
	\mathcal{H}_h=\int\frac{\rho_h\boldsymbol{u}_h\cdot\boldsymbol{u}_h}{2} + \rho_hgz + \frac{c_v}{c_p}\Theta_h\Pi_h\mathrm{d}\Omega\ .
\end{equation}

The diagnostic equations are spatially discretised in a similar manner by introducing test functions
$\boldsymbol{\beta}_h\in\mathbb{W}_1$, $\boldsymbol{w}_h\in\mathbb{W}_2$, $\xi_h,\sigma_h\in\mathbb{W}_3$ as
\begin{subequations}
	\begin{align}
		\int\overline{\rho}_h\boldsymbol{\beta}_h\cdot\overline{\boldsymbol{q}}_h\mathrm{d}\Omega =& \label{eq::diag_q}
		\int-(\nabla\times\boldsymbol{\beta}_h)\cdot\overline{\boldsymbol{u}}_h +
		\boldsymbol{\beta}_h\cdot f\hat{\boldsymbol{k}}\mathrm{d}\Omega
		-\int\boldsymbol{\beta}_h\cdot(\overline{\boldsymbol{u}}_h\cdot\hat{\boldsymbol{t}}_{\Lambda})\mathrm{d}\Lambda
		\\
		\int\boldsymbol{w}_h\cdot\overline{\boldsymbol{F}}_h\mathrm{d}\Omega =& \label{eq::diag_F}
		\int\frac{\boldsymbol{w}_h}{6}\cdot(2\rho_h^n\boldsymbol{u}_h^n +
		\rho_h^n\boldsymbol{u}_h^k + \rho_h^k\boldsymbol{u}_h^n + 2\rho_h^k\boldsymbol{u}_h^k)\mathrm{d}\Omega \\
		\int\xi_h\overline{\Phi}_h\mathrm{d}\Omega =& \int\frac{1}{6}\xi_h\Big(\boldsymbol{u}_h^n\cdot\boldsymbol{u}_h^n + \label{eq::diag_Phi}
		\boldsymbol{u}_h^n\cdot\boldsymbol{u}_h^k +
		\boldsymbol{u}_h^k\cdot\boldsymbol{u}_h^k\Big) + \xi_hgz\mathrm{d}\Omega \\
		\int\sigma_h\overline{\rho}_h\overline{\theta}_h\mathrm{d}\Omega =& \int\sigma_h\overline{\Theta}_h\mathrm{d}\Omega\ ,
	\end{align}
\end{subequations}
where $\Lambda$ is the domain outer boundary, and $\hat{\boldsymbol{t}}_{\Lambda}$ is the unit vector tangent to this boundary.

Second order accuracy in time is assured for \eqref{eq::ce_disc} via the exact second order integration of the variational 
derivatives of the energy, \eqref{eq::energy} with respect to the prognostic variables $\boldsymbol{u}_h$, $\rho_h$, $\Theta_h$ 
(respectively $\overline{\boldsymbol{F}}_h$, $\overline{\Phi}_h$, $\overline{\Pi}_h$), and the time centered representation of 
the terms within the nonlinear antisymmetric operators, $\overline{\boldsymbol{q}}_h$, $\overline{\theta}_h$ \cite{CH(2011)}.
This time centered discretisation results in an implicit, coupled system of nonlinear equations. Energy conservation is
derived by setting $\boldsymbol{v}_h=\overline{\boldsymbol{F}}_h$, $\phi_h=\overline{\Phi}_h$ and $\psi_h=\overline{\Pi}_h$
in \eqref{eq::ce_disc_u}, \eqref{eq::ce_disc_rho} and \eqref{eq::ce_disc_Theta} respectively. It is for this reason that we require the
collocation of the Exner pressure and density weighted potential temperature in order to ensure energy conservation.

While $\overline{\Pi}_h=(\Pi_h^n + \Pi_h^k)/2$ may be integrated exactly in time as a simple centered averaging for 
second order time integration, the Exner pressure itself is a non-polynomial function of the density weighted potential 
temperature \eqref{eq::Pi_cont}, which cannot be integrated exactly. Consequently we anticipate exact 
temporal conservation of energy only for a lowest order, piecewise constant representation of $\overline{\Pi}_h\in\mathbb{W}_3$. 
Previous results have shown the temporal loss of energy conservation 
due to the non-polynomial form of the equation of state at lowest order to be at machine precision \cite{Lee(2021)}.

Note that the discretisation presented in \eqref{eq::ce_disc} is generalised to a three dimensional 
domain. However the examples presented in Section \ref{sec::results} are all for one or two dimensional domains. For the one
dimensional domain, we have that $\overline{\boldsymbol{q}}_h:= 0$ in \eqref{eq::ce_disc_u}, such that
\eqref{eq::diag_q} is unnecessary. For the two dimensional case we have that the potential vorticity is
effectively a scalar function, $\overline{q}_h\in\mathbb{W}_0$, which is diagnosed subject to the 
test function $\beta_h\in\mathbb{W}_0$ in \eqref{eq::diag_q}. In all cases presented here there is
also no Coriolis term, such that $f:=0$.

\section{Helmholtz preconditioning}
\label{sec::helmholtz}

In this section we introduce a new preconditioner for the flux form evolution of the density weighted potential 
temperature on the Lorenz grid, as described in \eqref{eq::ce_disc}. 
This preconditioner is derived via the introduction of an auxiliary residual equation for the thermodynamic entropy
as a function of the density and density weighted potential temperature residuals. This allows us to express the Jacobian
operator in terms of the thermodynamic entropy instead of the density weighted potential temperature, which allows for
a linearisation of the coupled system that better captures the timescales of the buoyancy modes for a stratified atmosphere.

Following this derivation, we also describe the
other existing preconditioners which this will be compared against in Section \ref{sec::results}.
These include a previously derived Helmholtz operator for the flux form transport of $\Theta$ on the Lorenz grid 
\cite{Lee(2021)}
as described in Section \ref{sec::flux_orig_prec}, the material transport of $\theta$ on the Charney-Phillips grid 
\cite{Melvin(2019)},\cite{Maynard(2020)}
as given in Section \ref{sec::mat_prec}, and an energy conserving formulation of material form transport of $\theta$ on the 
Lorenz grid in Section \ref{sec::mat_prec_lorenz}.

\subsection{A novel flux form $\Theta$ preconditioner for the Lorenz grid}\label{sec::flux_new_prec}
In order to introduce the block preconditioner that we propose in this work, we first recall that in the
continuous form the thermodynamic entropy, $c_p\log(\theta)$ is subject to a material form advection equation \cite{LP(2021)},
and that this may be derived from the density and density weighted potential temperature flux 
form transport equations. Scaling \eqref{eq::rho_cont} and \eqref{eq::Theta_cont} respectively 
by $\rho^{-1}$ and $\Theta^{-1}$ and taking the difference of these expressions gives an expression
for the thermodynamic entropy transport as
\begin{subequations}\label{eq::eta_adv_cont}
	\begin{align}
	&\frac{1}{\Theta}\Big(\frac{\partial\Theta}{\partial t} + \nabla\cdot(\theta\boldsymbol{F})\Big) -
        \frac{1}{\rho}\Big(\frac{\partial\rho}{\partial t} + \nabla\cdot\boldsymbol{F}\Big) \\
        &= \frac{1}{\Theta}\frac{\partial\Theta}{\partial t} - \frac{1}{\rho}\frac{\partial\rho}{\partial t} +
        \frac{\boldsymbol{u}}{\theta}\cdot\nabla\theta  \\
	&= \frac{\partial}{\partial t}\log\Big(\frac{\Theta}{\rho}\Big) + \boldsymbol{u}\cdot\nabla\log(\theta) \\
	&= \frac{\partial\eta}{\partial t} + \boldsymbol{u}\cdot\nabla\eta = 0\ ,\label{eq::eta_adv_cont_4}
	\end{align}
\end{subequations}
where we define $\eta=\log(\theta)$ as the thermodynamic entropy scaled by $c_p^{-1}$.
We may therefore derive a discrete form of the thermodynamic entropy advection equation residual, $\mathcal{R}_{\eta}$
by first introducing the bilinear operators
\begin{subequations}
	\begin{align}
		\boldsymbol{\mathsf{M}}_3 &= \int\phi_h\psi_h\mathrm{d}\Omega\qquad\forall \phi_h,\psi_h\in\mathbb{W}_3 \\
		\boldsymbol{\mathsf{M}}_{3\rho} &= \int\rho_h\phi_h\psi_h\mathrm{d}\Omega\qquad\forall \phi_h,\psi_h\in\mathbb{W}_3 \\
		\boldsymbol{\mathsf{M}}_{3\Theta} &= \int\Theta_h\phi_h\psi_h\mathrm{d}\Omega\qquad\forall \phi_h,\psi_h\in\mathbb{W}_3\ .
	\end{align}
\end{subequations}
We then construct the residual error for the discrete thermodynamic entropy transport equation (with $\eta_h\in\mathbb{W}_3$) as
\begin{equation}\label{eq::R_eta}
	\mathcal{R}_{\eta} = \boldsymbol{\mathsf{M}}_3\Big(\boldsymbol{\mathsf{M}}_{3\Theta}^{-1}\mathcal{R}_{\Theta} - 
	\boldsymbol{\mathsf{M}}_{3\rho}^{-1}\mathcal{R}_{\rho}\Big)\ .
\end{equation}
Note that as the residual expressions $\mathcal{R}_{\Theta}$ and $\mathcal{R}_{\rho}$ converge, so too will the expression
for $\mathcal{R}_{\eta}$, provided that everywhere $\Theta_h$ and $\rho_h$ do not also converge to zero at a faster rate than
their corresponding residual expressions.
Note also that since basis functions spanning $\mathbb{W}_3$ are discontinuous
in all dimensions, the matrices $\boldsymbol{\mathsf{M}}_3$, $\boldsymbol{\mathsf{M}}_{3\rho}$ $\boldsymbol{\mathsf{M}}_{3\Theta}$
are block diagonal, and so their inverses may be computed directly with little expense.

We can now solve a coupled quasi-Newton problem for the increments at nonlinear iteration $k$: 
$\delta\boldsymbol{u}_h$, $\delta\rho_h$, $\delta\eta_h$, $\delta\Pi_h$ by an approximation
of the Jacobian as
\begin{equation}\label{eq::J}
	\begin{bmatrix}
		\boldsymbol{\mathsf{M}}_{2,R} & \boldsymbol{\mathsf{0}} & \boldsymbol{\mathsf{G}}_{\eta} & \boldsymbol{\mathsf{G}}_{\Pi} \\
		\boldsymbol{\mathsf{D}}_u & \boldsymbol{\mathsf{M}}_3 & \boldsymbol{\mathsf{0}} & \boldsymbol{\mathsf{0}} \\
		\boldsymbol{\mathsf{A}}_u & \boldsymbol{\mathsf{0}} & \boldsymbol{\mathsf{M}}_3 &  \boldsymbol{\mathsf{0}} \\
		\boldsymbol{\mathsf{0}} & \boldsymbol{\mathsf{C}}_{\rho} & \boldsymbol{\mathsf{C}}_{\eta} &  \boldsymbol{\mathsf{C}}_{\Pi} 
	\end{bmatrix}
	\begin{bmatrix}
		\delta\boldsymbol{u}_h \\ \delta\rho_h \\ \delta\eta_h \\ \delta\Pi_h
	\end{bmatrix} = -
	\begin{bmatrix}
		\mathcal{R}_u \\ \mathcal{R}_{\rho} \\ \mathcal{R}_{\eta} \\ \mathcal{R}_{\Pi}
	\end{bmatrix}\ ,
\end{equation}
for which the operators for the velocity equation Jacobian terms are given as
\begin{subequations}\label{eq::block_ops}
	\begin{align}
		\boldsymbol{\mathsf{M}}_2 =& \int\boldsymbol{w}_h\cdot\boldsymbol{v}_h\mathrm{d}\Omega
		\quad\forall\boldsymbol{w}_h,\boldsymbol{v}_h\in\mathbb{W}_2 \label{eq::M2}\\
		\boldsymbol{\mathsf{R}}_2 =& \frac{\Delta t}{2}\int\boldsymbol{w}_h\times f\boldsymbol{v}_h\mathrm{d}\Omega
		\quad\forall\boldsymbol{w}_h,\boldsymbol{v}_h\in\mathbb{W}_2 \label{eq::R2}\\
		\boldsymbol{\mathsf{G}}_{\eta} =& \frac{\Delta t}{2}\int\boldsymbol{v}_h\cdot\tilde{\nabla\Pi_h^n}\theta_h^n\phi_h\mathrm{d}\Omega 
		\label{eq::G_eta}
		\quad\forall\boldsymbol{v}_h\in\mathbb{W}_2,\phi_h\in\mathbb{W}_3 \\
		\boldsymbol{\mathsf{G}}_{\Pi} =& -\frac{\Delta t}{2}\int\nabla\cdot\boldsymbol{v}_h\theta_h^n\phi_h + 
		(\boldsymbol{v}_h\cdot\nabla\theta_h^n)\phi_h\mathrm{d}\Omega
		+\frac{\Delta t}{2}\int[\![\theta_h^n\boldsymbol{v}_h\cdot\hat{\boldsymbol{n}}_{\Gamma}]\!]\ldblbrace\phi_h\rdblbrace\mathrm{d}\Gamma
		\quad\forall\boldsymbol{v}_h\in\mathbb{W}_2,\phi_h\in\mathbb{W}_3\ ,
	\end{align}
\end{subequations}
and $\overset{\circ}{\boldsymbol{\mathsf{M}}}_2^{-1}$ is a lumped approximate inverse for the $\mathbb{W}_2$ mass matrix \eqref{eq::M2},
such that $\tilde{\nabla\Pi_h^n} = -\overset{\circ}{\boldsymbol{\mathsf{M}}}_2^{-1}\boldsymbol{\mathsf{D}}^{\top}\Pi_h^n$
is an approximate gradient of the Exner pressure. We also give the sum of the velocity space mass matrix and Coriolis 
operator as $\boldsymbol{\mathsf{M}}_{2,R} = \boldsymbol{\mathsf{M}}_2 + \boldsymbol{\mathsf{R}}_2$.
The operators for the density and entropy transport equation Jacobian terms are given as
\begin{subequations}\label{eq::block_ops_trans}
	\begin{align}
		\boldsymbol{\mathsf{D}} =& \int\phi_h\nabla\cdot\boldsymbol{v}_h\mathrm{d}\Omega
		\quad\forall\boldsymbol{v}_h\in\mathbb{W}_2,\phi_h\in\mathbb{W}_3 \\
		\boldsymbol{\mathsf{M}}_{2\rho} =& \int\rho_h^n\boldsymbol{w}_h\cdot\boldsymbol{v}_h\mathrm{d}\Omega
		\quad\forall\boldsymbol{w}_h,\boldsymbol{v}_h\in\mathbb{W}_2 \\
		\boldsymbol{\mathsf{D}}_u =& 
		\frac{\Delta t}{2}\boldsymbol{\mathsf{D}}\overset{\circ}{\boldsymbol{\mathsf{M}}}_2^{-1}\boldsymbol{\mathsf{M}}_{2\rho}\\
		\boldsymbol{\mathsf{A}}_u =& -\frac{\Delta t}{2}\int\nabla\cdot(\phi_h\boldsymbol{v}_h)\eta_h^n\mathrm{d}\Omega \label{eq::A_u}
		+\frac{\Delta t}{2}\int[\![\phi_h\boldsymbol{v}_h\cdot\hat{\boldsymbol{n}}_{\Gamma}]\!]\ldblbrace\eta_h^n\rdblbrace\mathrm{d}\Gamma
		\quad\forall\boldsymbol{v}_h\in\mathbb{W}_2,\phi_h\in\mathbb{W}_3\,
	\end{align}
\end{subequations}
while those for the equation of state are given as
\begin{subequations}\label{eq::block_ops_eos}
	\begin{align}
		\boldsymbol{\mathsf{C}}_{\rho} =& 
		-\frac{R}{c_v}\boldsymbol{\mathsf{M}}_3\boldsymbol{\mathsf{M}}_{3\rho}^{-1}\boldsymbol{\mathsf{M}}_3 \label{eq::C_rho}\\
		\boldsymbol{\mathsf{C}}_{\eta} =&  
		-\frac{R}{c_v}\boldsymbol{\mathsf{M}}_3 \label{eq::C_eta}\\
		\boldsymbol{\mathsf{M}}_{3\Pi} =& \int\Pi_h\phi_h\psi_h\mathrm{d}\Omega\qquad\forall \phi_h,\psi_h\in\mathbb{W}_3 \\
		\boldsymbol{\mathsf{C}}_{\Pi} =& 
		\boldsymbol{\mathsf{M}}_3\boldsymbol{\mathsf{M}}_{3\Pi}^{-1}\boldsymbol{\mathsf{M}}_3.
	\end{align}
\end{subequations}
Since the equations of motion are a function of $\Theta_h$ and not $\eta_h$ \eqref{eq::ce_disc}, transformations 
are required in order to derive the operators $\boldsymbol{\mathsf{G}}_{\eta}$, 
$\boldsymbol{\mathsf{C}}_{\eta}$ and $\boldsymbol{\mathsf{C}}_{\eta}$ 
in \eqref{eq::G_eta}
\eqref{eq::C_rho} and \eqref{eq::C_eta} respectively. 
The operator $\boldsymbol{\mathsf{G}}_{\eta}$ is derived by taking the variational derivative of the term 
$\frac{\delta}{\delta\eta_h}\int\boldsymbol{v}_he^{\eta_h^n}\tilde{\nabla\Pi_h^n}\mathrm{d}\Omega=
\int\boldsymbol{v}_h\theta_h^n\tilde{\nabla\Pi_h^n}\phi_h\mathrm{d}\Omega$
$\forall\boldsymbol{v}_h\in\mathbb{W}_2,\phi_h\in\mathbb{W}_3$. The operators 
$\boldsymbol{\mathsf{C}}_{\eta}$ and $\boldsymbol{\mathsf{C}}_{\eta}$ are derived by expanding the term in the
equation of state for $\Theta_h$ \eqref{eq::ce_disc_Pi} as
$\int\chi_h\frac{R}{c_v}\log(\Theta_h)\mathrm{d}\Omega=
\int\chi_h\frac{R}{c_v}(\log(\rho_h) + \log(\theta_h))\mathrm{d}\Omega=
\int\chi_h\frac{R}{c_v}(\log(\rho_h) + \eta_h)\mathrm{d}\Omega$. The resulting operators are then derived by 
differentiating this with respect to $\rho_h$ and $\eta_h$.

The approximate Jacobian in \eqref{eq::J} has the same block structure as the operator for the {\tt LFRic} model 
\cite{Melvin(2019)},\cite{Maynard(2020)}, which solves for the material transport of potential temperature on the Charney-Phillips grid
(the vertical component of the $\mathbb{W}_2$ space), not the flux form transport of density weighted potential temperature
on the Lorenz grid as is done above. However the precise form of the blocks $\boldsymbol{\mathsf{G}}_{\eta}$ and 
$\boldsymbol{\mathsf{C}}_{\eta}$ are different from their equivalents in {\tt LFRic} (as is the space for the $\delta\eta_h$ 
test functions).

Like the {\tt LFRic} preconditioner, successive block factorisation is applied as
\begin{subequations}\label{eq::incs}
	\begin{align}
		\delta\eta_h &= -\boldsymbol{\mathsf{M}}_3^{-1}(\mathcal{R}_{\eta} + \boldsymbol{\mathsf{A}}_u\delta\boldsymbol{u}_h) \\
		\delta\rho_h &= -\boldsymbol{\mathsf{M}}_3^{-1}(\mathcal{R}_{\rho} + \boldsymbol{\mathsf{D}}_u\delta\boldsymbol{u}_h) \\
		\delta\boldsymbol{u}_h &= -\overset{\circ}{\boldsymbol{\mathcal{M}}}_2^{-1}(\mathcal{R}_u' + \boldsymbol{\mathsf{G}}_{\Pi}\delta\Pi_h)\ ,
	\end{align}
\end{subequations}
where $\boldsymbol{\mathcal{M}}_2=\boldsymbol{\mathsf{M}}_{2,R} - 
\boldsymbol{\mathsf{G}}_{\eta}\boldsymbol{\mathsf{M}}_3^{-1}\boldsymbol{\mathsf{A}}_u$, 
$\overset{\circ}{\boldsymbol{\mathcal{M}}}_2^{-1}$ is its lumped approximate inverse, and
$\mathcal{R}_u' = \mathcal{R}_u - \boldsymbol{\mathsf{G}}_{\eta}\boldsymbol{\mathsf{M}}_3^{-1}\mathcal{R}_{\eta}$. This results 
in a Helmholtz equation for the Exner pressure increment as
\begin{equation}\label{eq::H}
	\begin{bmatrix}
		\boldsymbol{\mathsf{C}}_{\Pi} + 
		(\boldsymbol{\mathsf{C}}_{\rho}\boldsymbol{\mathsf{M}}_3^{-1}\boldsymbol{\mathsf{D}}_u - \frac{R}{c_v}\boldsymbol{\mathsf{A}}_u)
		\overset{\circ}{\boldsymbol{\mathcal{M}}}_2^{-1}\boldsymbol{\mathsf{G}}_{\Pi}
	\end{bmatrix}\delta\Pi_h = 
		-\mathcal{R}_{\Pi} + 
		\boldsymbol{\mathsf{C}}_{\rho}\boldsymbol{\mathsf{M}}_3^{-1}\mathcal{R}_{\rho} -
		\frac{R}{c_v}\mathcal{R}_{\eta} 
		-(\boldsymbol{\mathsf{C}}_{\rho}\boldsymbol{\mathsf{M}}_3^{-1}\boldsymbol{\mathsf{D}}_u - \frac{R}{c_v}\boldsymbol{\mathsf{A}}_u)
		\overset{\circ}{\boldsymbol{\mathcal{M}}}_2^{-1}\mathcal{R}_u'.
\end{equation}
Once the Exner pressure increment at the current Newton iteration has been determined, the increments for the other variables,
$\delta\eta_h$, $\delta\rho_h$, $\delta\boldsymbol{u}_h$ can be determined respectively via \eqref{eq::incs}. Recalling that our
actual prognostic variable is the density weighted potential temperature, $\Theta_h^k$, and not the thermodynamic entropy, $\eta_h^k$, 
this may then be determined at Newton iteration $k$ $\forall\sigma_h,\psi_h,\phi_h\in\mathbb{W}_3$ as
\begin{subequations}
\begin{align}
	\int\sigma_h\rho_h^{k-1}\theta_h^{k-1}\mathrm{d}\Omega &= \int\sigma_h\Theta_h^{k-1}\mathrm{d}\Omega,\\
	\int\psi_h\eta_h^k\mathrm{d}\Omega &= \int\psi_h(\log(\theta_h^{k-1}) + \delta\eta_h^k)\mathrm{d}\Omega,\label{eq::eta_disc}\\
	\int\phi_h\Theta_h^k\mathrm{d}\Omega &= \int\phi_h\rho_h^ke^{\eta_h^k}\mathrm{d}\Omega\ .
\end{align}
\end{subequations}
Note that since $\mathcal{R}_{\eta}$ is an auxiliary expression only \eqref{eq::R_eta}, we only ever require $\eta_h^k$ 
as derived from \eqref{eq::eta_disc} in order to update $\Theta_h^k$. However this is not required in order to evaluate the original 
equations of motion \eqref{eq::ce_disc}.
Also note that as stated previously no Coriolis term is applied in the example configurations presented below, such that the
operator $\boldsymbol{\mathsf{R}}_2$ \eqref{eq::R2} is omitted.

In this study we only ever solve for the Helmholtz equation, as given in \eqref{eq::H} and the analogous forms below,
and never for the full coupled system as given in \eqref{eq::J}. However for some applications the mass lumping in \eqref{eq::H}
may prove ineffective, and so an outer solve of the coupled system may be necessary. 
In such a care the Helmholtz solver may still be used as a preconditioner to accelerate the solution of the outer solve.
Particular examples where lumping of the velocity mass matrix inverse may lead to poor approximations of its inverse
and degraded accuracy for the resulting Helmholtz operator may include distorted geometries with large non-orthogonal 
vector components, and higher order bases with less compact mass matrix operators.

Similar to the {\tt LFRic} preconditioner \cite{Maynard(2020)} discussed in Section \ref{sec::mat_prec}, the new flux form preconditioner also scales 
temporally with the acoustic and buoyancy modes. Inspecting the block matrices in \eqref{eq::block_ops} we see that the 
Helmholtz operator in \eqref{eq::H} scales as
\begin{equation}\label{eq::H_scaling_1}
	\frac{1}{\Pi_h^n} + 
	\frac{R}{c_v}\Bigg(\frac{\Delta t}{\rho_h^n}\frac{\partial\rho_h^n}{\partial x} + \Delta t\frac{\partial\eta_h^n}{\partial x}\Bigg)
	\Bigg(1 + \Delta t^2\theta_h^n\frac{\partial\Pi_h^n}{\partial x}\frac{\partial\eta_h^n}{\partial x}\Bigg)^{-1}
	\frac{\Delta t}{\Delta x}\theta_h^n\ .
\end{equation}
Recalling the equation of state \eqref{eq::Pi_cont} and that the pressure is determined from the Exner pressure
as $p = p_0(\Pi/c_p)^{c_p/R}$ we have the square of the speed of sound as $c_s^2=\delta p/\delta\rho=R\theta\Pi/c_v$.
From the hydrostatic balance relation we also have an approximate equality as $g=-\theta\partial\Pi/\partial x$.
We also have the square of the Brunt-V\"ais\"al\"a frequency as 
$\mathcal{N}^2=(g/\theta)(\partial\theta/\partial x)=g\partial\eta/\partial x$. Combining these relations and scaling
\eqref{eq::H_scaling_1} by $\Pi_h^n$ gives
\begin{equation}\label{eq::H_scaling_2}
	1 + c_s^2
	\Bigg(\frac{\Delta t}{\Delta x}\Bigg)^2
	\Bigg(\frac{\Delta x}{\rho_h^n}\frac{\partial\rho_h^n}{\partial x} + \frac{\Delta x}{g}\mathcal{N}^2\Bigg)
	\Bigg(1 - \Delta t^2\mathcal{N}^2\Bigg)^{-1}\ .
\end{equation}
This scaling of the Helmholtz operator with both $c_s^2$ and $\mathcal{N}^2$ is also seen for the {\tt LFRic} 
preconditioner \cite{Maynard(2020)}, however in the present case this is achieved via the relation 
$\partial\eta/\partial x=\theta^{-1}(\partial\theta/\partial x)$. 
Since the scaling with respect to $\mathcal{N}^2$ arises due to the operator
$\boldsymbol{\mathsf{A}}_u$ \eqref{eq::A_u}, it may be possible to simplify the structure of the Helmholtz
operator somewhat without degrading the performance significantly by only accounting for gradients and
boundary integrals in the vertical for this term. The operator also scales with spatial variations in 
density as $\partial\log(\rho_h^n)/\partial x$. We also note that the numerator $1-\Delta t^2\mathcal{N}^2$ 
is lumped so that the spatial derivative will not be properly approximated here. For time steps of more than 
a couple of minutes with a standard atmospheric stratification, as described in the Appendix and \cite{UMJS(2014)}, 
we also have $\Delta t^2\mathcal{N}^2>>1$.

\subsection{Original flux form $\Theta$ preconditioner for the Lorenz grid}\label{sec::flux_orig_prec}
In a previous work, an alternative preconditioner was presented for the flux form transport of $\Theta$
on the Lorenz grid \cite{Lee(2021)}, and implemented in the vertical dimension as part of a horizontally 
explicit, vertically implicit scheme for the 3D compressible Euler equations using compatible finite elements.
While the precise form of the operators is detailed in \cite{Lee(2021)}, the block structure of this preconditioner
is given as
\begin{equation}\label{eq::J_lee21}
	\begin{bmatrix}
		\boldsymbol{\mathsf{M}}_{2,R} & \boldsymbol{\mathsf{0}} & \boldsymbol{\mathsf{G}}_{\Theta} & \boldsymbol{\mathsf{G}}_{\Pi} \\
		\boldsymbol{\mathsf{D}}_u & \boldsymbol{\mathsf{M}}_3 & \boldsymbol{\mathsf{0}} & \boldsymbol{\mathsf{0}} \\
		\boldsymbol{\mathsf{D}}_{\Theta} & \boldsymbol{\mathsf{Q}}_u & \boldsymbol{\mathsf{M}}_3 &  \boldsymbol{\mathsf{0}} \\
		\boldsymbol{\mathsf{0}} & \boldsymbol{\mathsf{0}} & \boldsymbol{\mathsf{C}}_{\Theta} &  \boldsymbol{\mathsf{C}}_{\Pi} 
	\end{bmatrix}
	\begin{bmatrix}
		\delta\boldsymbol{u}_h \\ \delta\rho_h \\ \delta\Theta_h \\ \delta\Pi_h
	\end{bmatrix} = -
	\begin{bmatrix}
		\mathcal{R}_u \\ \mathcal{R}_{\rho} \\ \mathcal{R}_{\Theta} \\ \mathcal{R}_{\Pi}
	\end{bmatrix}\ .
\end{equation}
Note that in contrast to the preconditioner used to solve for $\delta\eta_h$ \eqref{eq::J}, the above operator has 
no term in the $[\Pi,\rho]$ block, and instead has a non-zero term in the $[\Theta,\rho]$ block, 
$\boldsymbol{\mathsf{Q}}_u$. Repeated Schur complement reduction then leads to an operator of the form
\begin{multline}\label{eq::H2}
	\big[
		\boldsymbol{\mathsf{M}}_3 - 
		(\boldsymbol{\mathsf{D}}_{\Theta} - \boldsymbol{\mathsf{Q}}_{\Theta}\boldsymbol{\mathsf{M}}_3^{-1}\boldsymbol{\mathsf{D}}_{\rho})
		\overset{\circ}{\boldsymbol{\mathsf{M}}}_{2,R}^{-1} 
		(\boldsymbol{\mathsf{G}}_{\Theta} - \boldsymbol{\mathsf{G}}_{\Pi}\boldsymbol{\mathsf{C}}_{\Pi}^{-1}\boldsymbol{\mathsf{C}}_{\Theta})
	\big]\delta\Theta_h = \\
	-\mathcal{R}_{\Theta} + \boldsymbol{\mathsf{Q}}_{\Theta}\boldsymbol{\mathsf{M}}_3^{-1}\mathcal{R}_{\rho} +
	(\boldsymbol{\mathsf{D}}_{\Theta} -\boldsymbol{\mathsf{Q}}_{\Theta}\boldsymbol{\mathsf{M}}_3^{-1}\boldsymbol{\mathsf{D}}_{\rho})
		\overset{\circ}{\boldsymbol{\mathsf{M}}}_{2,R}^{-1}(\mathcal{R}_u - \boldsymbol{\mathsf{G}}_{\Pi}\boldsymbol{\mathsf{C}}_{\Pi}^{-1}\mathcal{R}_{\Pi})\ .
\end{multline}
The benefit of \eqref{eq::H2} with respect to \eqref{eq::H} is that the structure of the lumped inverse blocks is simpler, such that
more of the dynamics is accounted for within the Helmholtz operator itself. However the downside is that the block $\boldsymbol{\mathsf{Q}}_u$
scales with the velocity, $\boldsymbol{u}$, as opposed to the block $\boldsymbol{\mathsf{C}}_{\rho}$ in \eqref{eq::H}, which scales with 
the inverse density. Since in most applications the vertical motions of the atmosphere are small, this block will have a negligible 
contribution to the solution in many cases.

\subsection{Material form $\theta$ preconditioner for the Charney-Phillips grid ({\tt LFRic})}\label{sec::mat_prec}

The preconditioner in \eqref{eq::H} has a similar structure as one for the material form transport of $\theta$
on the Charney-Phillips grid \cite{Melvin(2019)},\cite{Maynard(2020)}. In that case the approximate Jacobian is given 
as
\begin{equation}\label{eq::J_LFRic}
	\begin{bmatrix}
		\boldsymbol{\mathsf{M}}_{2,R} & \boldsymbol{\mathsf{0}} & \boldsymbol{\mathsf{G}}_{\theta} & \boldsymbol{\mathsf{G}}_{\Pi} \\
		\boldsymbol{\mathsf{D}}_u & \boldsymbol{\mathsf{M}}_3 & \boldsymbol{\mathsf{0}} & \boldsymbol{\mathsf{0}} \\
		\boldsymbol{\mathsf{A}}_{\theta,u} & \boldsymbol{\mathsf{0}} & \boldsymbol{\mathsf{M}}_{\theta} &  \boldsymbol{\mathsf{0}} \\
		\boldsymbol{\mathsf{0}} & \boldsymbol{\mathsf{C}}_{\rho} & \boldsymbol{\mathsf{C}}_{\theta} &  \boldsymbol{\mathsf{C}}_{\Pi} 
	\end{bmatrix}
	\begin{bmatrix}
		\delta\boldsymbol{u}_h \\ \delta\rho_h \\ \delta\theta_h \\ \delta\Pi_h
	\end{bmatrix} = -
	\begin{bmatrix}
		\mathcal{R}_u \\ \mathcal{R}_{\rho} \\ \mathcal{R}_{\theta} \\ \mathcal{R}_{\Pi}
	\end{bmatrix},
\end{equation}
which reduces to a Helmholtz equation of the form
\begin{multline}\label{eq::H_LFRic}
	\begin{bmatrix}
		\boldsymbol{\mathsf{C}}_{\Pi} + 
		(\boldsymbol{\mathsf{C}}_{\theta}\overset{\circ}{\boldsymbol{\mathsf{M}}}_{\theta}^{-1}\boldsymbol{\mathsf{A}}_{\theta,u} + 
		\boldsymbol{\mathsf{C}}_{\rho}\boldsymbol{\mathsf{M}}_3^{-1}\boldsymbol{\mathsf{D}}_u)
		\overset{\circ}{\boldsymbol{\mathcal{M}}}_{2,\theta}^{-1}\boldsymbol{\mathsf{G}}_{\Pi}
	\end{bmatrix}\delta\Pi_h = \\
		-\mathcal{R}_{\Pi} + 
		\boldsymbol{\mathsf{C}}_{\theta}\overset{\circ}{\boldsymbol{\mathsf{M}}}_{\theta}^{-1}\mathcal{R}_{\theta} + 
		\boldsymbol{\mathsf{C}}_{\rho}\boldsymbol{\mathsf{M}}_3^{-1}\mathcal{R}_{\rho} -
		(\boldsymbol{\mathsf{C}}_{\theta}\overset{\circ}{\boldsymbol{\mathsf{M}}}_{\theta}^{-1}\boldsymbol{\mathsf{A}}_{\theta,u} + 
		\boldsymbol{\mathsf{C}}_{\rho}\boldsymbol{\mathsf{M}}_3^{-1}\boldsymbol{\mathsf{D}}_u)
		\overset{\circ}{\boldsymbol{\mathcal{M}}}_{2,\theta}^{-1}\mathcal{R}_u'\ ,
\end{multline}
where $\boldsymbol{\mathcal{M}}_{2,\theta}=\boldsymbol{\mathsf{M}}_{2,R} - 
\boldsymbol{\mathsf{G}}_{\theta}\overset{\circ}{\boldsymbol{\mathsf{M}}}_{\theta}^{-1}\boldsymbol{\mathsf{A}}_{\theta,u}$.

One of the benefits of flux form transport of $\Theta_h$ is that the potential temperature variance, 
$\mathcal{Z}_h=\int\Theta_h^2/(2\rho_h)\mathrm{d}\Omega$
is a mathematical entropy of the dry compressible Euler equations (all eigenvalues $\ge 0$), and so conserving or provably 
damping $\mathcal{Z}_h$ can help to stabilise thermal processes \cite{Ricardo(2024)}.
Note that this mathematical entropy is distinct from the thermodynamic entropy that we have previously discussed.
Therefore it is tempting to look at solving for flux form $\Theta$ on the Charney-Phillips grid, in order to stabilise 
thermal processes without incurring the spurious computational mode that is associated with the Lorenz staggering 
\cite{Arakawa(1996)}. However doing so leads to the inconsistent material transport of $\rho_h\in\mathbb{W}_3$ 
\cite{Eldred(2019)}, resulting in convergence problems for the nonlinear solver. 
As discussed in the introduction, this issue may potentially be 
negated however via the rehabilitation of the horizontal fluxes for the density weighted potential temperature equation, 
as has previously been applied in order to recover conservation on the Charney-Philips grid using finite differences
\cite{Thuburn(2022)}.

\subsection{Material form $\theta$ preconditioner for the Lorenz grid}\label{sec::mat_prec_lorenz}

One disadvantage of the Charney-Phillips grid is that by staggering the Exner pressure and potential 
temperature energy conservation is not preserved discretely. This is because the variational derivative of the
energy with respect to the potential temperature is a function of the Exner pressure, and so 
representing these on different spaces breaks the anti-symmetry of the Hamiltonian structure of the 
energy conserving formulation \cite{Eldred(2019)},\cite{Lee(2021)}. 

Comparing preconditioners for 
energy conserving and non-conserving spatial discretisations is somewhat problematic, since the 
energy conserving formulation is innately more stable, separate from the choice of linearisations 
made in the approximation of the Jacobian operator. For a more thorough comparison we therefore also 
introduce an energy conserving formulation with material form transport of $\theta_h$ collocated 
with $\Pi_h\in\mathbb{W}_3$ (the Lorenz grid), for which the approximate Jacobian has the same form
as for the material form transport of potential temperature on the Charney-Phillips grid 
\cite{Melvin(2019)},\cite{Maynard(2020)}, albeit with a different finite element space for the 
potential temperature. The energy conserving formulation for material form potential temperature 
transport on the Lorenz grid is given 
for $\boldsymbol{v}_h\in\mathbb{W}_2$, $\phi_h,\psi_h,\chi_h\in\mathbb{W}_3$ as

\begin{subequations}\label{eq::ce_disc_matadv_lorenz}
\begin{align}
	\int\boldsymbol{v}_h\cdot(\boldsymbol{u_h}^k-\boldsymbol{u}_h^n) + 
	\Delta t\boldsymbol{v}_h\cdot\overline{\boldsymbol{q}}_h\times\overline{\boldsymbol{F}}_h -
	\Delta t(\nabla\cdot\boldsymbol{v}_h)\overline{\Phi}_h^L +
	\Delta t\nabla\cdot\Big(\frac{\boldsymbol{v}_h\overline{\rho\Pi}_h^L}{\overline{\rho}_h}\Big)\overline{\theta}_h\mathrm{d}\Omega \\\notag
	-\Delta t\int
	[\![\frac{\overline{\rho\Pi}_h^L\boldsymbol{v}_h\cdot\hat{\boldsymbol{n}}_{\Gamma}}{\overline{\rho}_h}]\!]\ldblbrace\overline{\theta}_h\rdblbrace 
	+ c
	[\![\frac{\overline{\rho\Pi}_h^L\boldsymbol{v}_h\cdot\hat{\boldsymbol{n}}_{\Gamma}}{\overline{\rho}_h}]\!][\![\overline{\theta}_h]\!] + 
	\mathrm{d}\Gamma
	&= \mathcal{R}_u^L\\
	\int\phi_h(\rho_h^k-\rho_h^n) + \Delta t\phi_h\nabla\cdot\overline{\boldsymbol{F}}_h\mathrm{d}\Omega& = \mathcal{R}_{\rho}^L \\
	\int\psi_h(\theta^k_h-\theta^n_h) -
	\Delta t\nabla\cdot\Big(\frac{\psi_h\overline{\boldsymbol{F}}_h}{\overline{\rho}_h}\Big)\overline{\theta}_h 
	\mathrm{d}\Omega +
	\Delta t\int
	[\![\frac{\psi_h\overline{\boldsymbol{F}}_h\cdot\hat{\boldsymbol{n}}_{\Gamma}}{\overline{\rho}_h}]\!]\ldblbrace\overline{\theta}_h\rdblbrace
	+ c
	[\![\frac{\psi_h\overline{\boldsymbol{F}}_h\cdot\hat{\boldsymbol{n}}_{\Gamma}}{\overline{\rho}_h}]\!][\![\overline{\theta}_h]\!]\mathrm{d}\Gamma 
	&= \mathcal{R}_{\theta}^L\\
	\int\chi_h\log(\Pi_h^k) - \chi_h\frac{R}{c_v}\log(\rho_h^k) -
	- \chi_h\frac{R}{c_v}\log(\theta_h^k) -
	\chi_h\frac{R}{c_v}\log\Big(\frac{R}{p_0}\Big) -
	\chi_h\log(c_p)\mathrm{d}\Omega &= \mathcal{R}_{\Pi}^L\ ,
\end{align}
\end{subequations}

\noindent
where $\overline{\Phi}_h^L$ and $\overline{\rho\Pi}_h^L$ are the variational derivatives of the energy 
$\mathcal{H}_h^L=\int\frac{1}{2}\rho_h\boldsymbol{u}_h\cdot\boldsymbol{u}_h + \rho_hgz + \frac{c_v}{c_p}\rho_h\theta_h\Pi_h\mathrm{d}\Omega$
with respect to $\rho_h$ and $\theta_h$. The above formulation has the same antisymmetric structure as 
for a previous energy conserving form of the thermal shallow water equations \cite{Eldred(2019)}, and 
conserves energy for a choice of $\boldsymbol{v}_h=\overline{\boldsymbol{F}}_h$, $\phi_h=\overline{\Phi}_h^L$,
$\psi_h=\overline{\rho\Pi}_h^L$.

The block structure of the material form $\theta$ transport Jacobian is given as
\begin{equation}
	\begin{bmatrix}
		\boldsymbol{\mathsf{M}}_{2,R} & \boldsymbol{\mathsf{0}} & \boldsymbol{\mathsf{G}}_{\theta} & \boldsymbol{\mathsf{G}}_{\Pi} \\
		\boldsymbol{\mathsf{D}}_u & \boldsymbol{\mathsf{M}}_3 & \boldsymbol{\mathsf{0}} & \boldsymbol{\mathsf{0}} \\
		\boldsymbol{\mathsf{A}}_u & \boldsymbol{\mathsf{0}} & \boldsymbol{\mathsf{M}}_3 &  \boldsymbol{\mathsf{0}} \\
		\boldsymbol{\mathsf{0}} & \boldsymbol{\mathsf{C}}_{\rho} & \boldsymbol{\mathsf{C}}_{\theta} &  \boldsymbol{\mathsf{C}}_{\Pi} 
	\end{bmatrix}
	\begin{bmatrix}
		\delta\boldsymbol{u}_h \\ \delta\rho_h \\ \delta\theta_h \\ \delta\Pi_h
	\end{bmatrix} = -
	\begin{bmatrix}
		\mathcal{R}_u \\ \mathcal{R}_{\rho} \\ \mathcal{R}_{\theta} \\ \mathcal{R}_{\Pi}
	\end{bmatrix}.
\end{equation}
This structure is similar to that for the new preconditioner described in Section 3.1 \eqref{eq::J},
with the principle differences being the operators $\boldsymbol{\mathsf{G}}_{\theta}$, $\boldsymbol{\mathsf{C}}_{\theta}$, which 
in contrast to the analogous operators for the material form $\eta$ system \eqref{eq::G_eta}, \eqref{eq::C_eta} are given here as
\begin{subequations}
	\begin{align}
		\boldsymbol{\mathsf{G}}_{\theta} =& \frac{\Delta t}{2}\int\boldsymbol{v}_h\cdot\tilde{\nabla\Pi_h^n}\theta_h^n\phi_h\mathrm{d}\Omega 
		\quad\forall\boldsymbol{v}_h\in\mathbb{W}_2,\phi_h\in\mathbb{W}_3 \\
		\boldsymbol{\mathsf{M}}_{3\theta} =& \int\theta_h\phi_h\psi_h\mathrm{d}\Omega\qquad\forall \phi_h,\psi_h\in\mathbb{W}_3 \\
		\boldsymbol{\mathsf{C}}_{\theta} =& -\frac{R}{c_v}
		\boldsymbol{\mathsf{M}}_3\boldsymbol{\mathsf{M}}_{3\theta}^{-1}\boldsymbol{\mathsf{M}}_3.
	\end{align}
\end{subequations}
This different block structure leads to a slightly different Helmholtz operator from that described in \eqref{eq::H} as
\begin{multline}
	\begin{bmatrix}
		\boldsymbol{\mathsf{C}}_{\Pi} + 
		(\boldsymbol{\mathsf{C}}_{\rho}\boldsymbol{\mathsf{M}}_3^{-1}\boldsymbol{\mathsf{D}}_u +
		\boldsymbol{\mathsf{C}}_{\theta}\boldsymbol{\mathsf{M}}_{3}^{-1}\boldsymbol{\mathsf{A}}_u)
		\overset{\circ}{\boldsymbol{\mathcal{M}}}_2^{-1}\boldsymbol{\mathsf{G}}_{\Pi}
	\end{bmatrix}\delta\Pi_h = \\
		-\mathcal{R}_{\Pi} + 
		\boldsymbol{\mathsf{C}}_{\rho}\boldsymbol{\mathsf{M}}_3^{-1}\mathcal{R}_{\rho} +
		\boldsymbol{\mathsf{C}}_{\theta}\boldsymbol{\mathsf{M}}_{3}^{-1}\mathcal{R}_{\theta}
		-(\boldsymbol{\mathsf{C}}_{\rho}\boldsymbol{\mathsf{M}}_3^{-1}\boldsymbol{\mathsf{D}}_u +
		\boldsymbol{\mathsf{C}}_{\theta}\boldsymbol{\mathsf{M}}_{3}^{-1}\boldsymbol{\mathsf{A}}_u)
		\overset{\circ}{\boldsymbol{\mathcal{M}}}_2^{-1}\mathcal{R}_u'.
\end{multline}

\section{Results}
\label{sec::results}

The numerical algorithms in Section~\ref{sec::mixed_fem} and~\ref{sec::helmholtz} were implemented in the Julia programming language. In particular, the mixed finite element discretisation of the compressible Euler equations 
and the computation of the different blocks in the approximate block Jacobian matrices were implemented 
using the \texttt{Gridap} \cite{Badia(2020)} finite element framework. 
This package provides a rich set of software tools in order to define and evaluate 
the cell and facet integrals terms in the different weak discrete variational formulations using a highly expressive and compact syntax that resembles mathematical notation. It also provides the tools required to assemble these terms 
selectively and flexibly into different sparse matrices and vectors as per required by the underlying block preconditioned iterative solvers. 

For the 2D experiments below, we also used the \texttt{GridapDistributed} package \cite{Badia(2022)} to implement the distributed-memory parallelization of the algorithms at hand. For efficiency and scalability in mind, the implemented message-passing code does not build explicitly the Helmholtz operator resulting from repeated Schur complement reduction (this would require, among others, sparse matrix-matrix multiplications, and extra memory consumption). Instead, it uses approximate lumping of the velocity mass matrix inverse, codes the action of the Helmholtz operator on a given vector (without building it explicitly in memory), and uses this action to solve the linear system using a preconditioned GMRES linear solver iteration. 

All numerical experiments in this section were conducted on the Gadi petascale supercomputer, hosted by the Australian National Computational Infrastructure (NCI). All floating-point operations are performed in IEEE double precision.

\subsection{1D atmosphere with potential temperature perturbation}

\begin{figure}[!hbtp]
\begin{center}
\includegraphics[width=0.48\textwidth,height=0.36\textwidth]{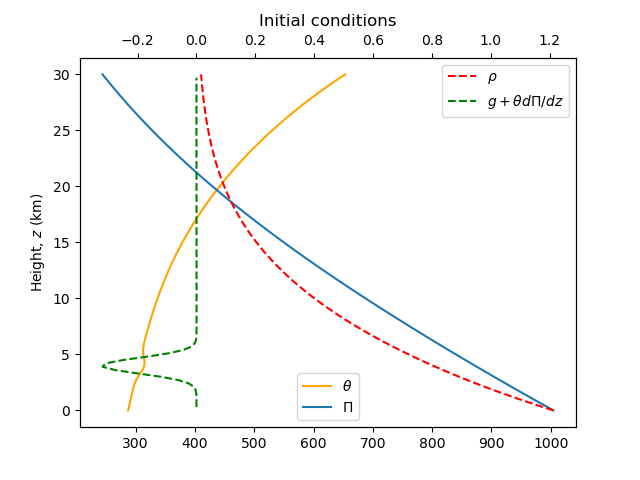}
	\caption{Initial vertical profiles for the 1D test case, with a potential 
	temperature perturbation overlaying a state of hydrostatic balance.}
\label{fig::vert_profs}
\end{center}
\end{figure}

The stability and efficiency of the different preconditioners detailed in
Section \ref{sec::helmholtz} were compared for the solution of a stratified 1D atmosphere with a height of $30km$, with vertical 
profiles similar to those used in an existing test case for baroclinic instability on the sphere \cite{UMJS(2014)}, 
as given in the Appendix. 
While in that test the initial vertical profiles are specified so as to establish a 
state of hydrostatic balance, here we perturb the potential temperature field by a Gaussian bubble of the form
\begin{equation}
	\theta_p = 10\exp(-10^{-6}(z-4000)^2).
\end{equation}
This perturbation has the effect of generating internal waves that radiate out from the bubble (which itself stays
relatively stationary due to the positive lapse rate of the initial conditions), so as to provide a challenging 
dynamical process for the different formulations of the preconditioner to accommodate. The initial conditions are given in Figure 
\ref{fig::vert_profs}, including the term $g + \theta d\Pi/dz$, which reflects the deviation from hydrostatic 
balance as generated by $\theta_p$.

In all cases we use a time step of $\Delta t=600s$ with 100 uniformly spaced lowest order elements, and
boundary conditions on the bottom $z_b$ and top $z_t$ of the form $w|_{z_b}=w|_{z_t}=0$, 
$\partial\Pi/\partial z|_{z_b} = \partial\Pi/\partial z|_{z_t}=0$.
No upwinding or other form of stabilisation was applied in any case, such that $c=0$ in 
\eqref{eq::ce_disc_u} and \eqref{eq::ce_disc_Theta}. 
The addition of upwinding or additional damping terms would serve to suppress
model instability in the cases where this is observed below. However our focus here is in the
inherent representation of dynamics and model stability for the formulations and preconditioners presented 
in Section \ref{sec::helmholtz}. In any case, the implementation of such additional 
terms would differ between the formulations, which would somewhat compromise the comparisons 
presented below.

We explore results
for two different formulations, one in which the matrix inverses are all computed directly, and the
solvers are iterated to convergence, and a second using lumped matrix inverses (with the exception
of the Helmholtz operator, which is solved directly), and just four Newton nonlinear iterations per time step, since in practice for 
high resolution simulations using parallel decompositions direct matrix inverses are not practical. In each case the 
Helmholtz operator and all other matrices are assembled only once at the beginning of each time step.

\begin{figure}[!hbtp]
\begin{center}
\includegraphics[width=0.32\textwidth,height=0.24\textwidth]{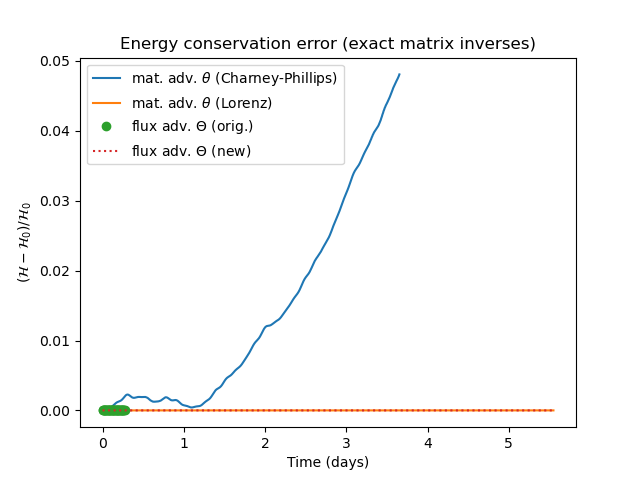}
\includegraphics[width=0.32\textwidth,height=0.24\textwidth]{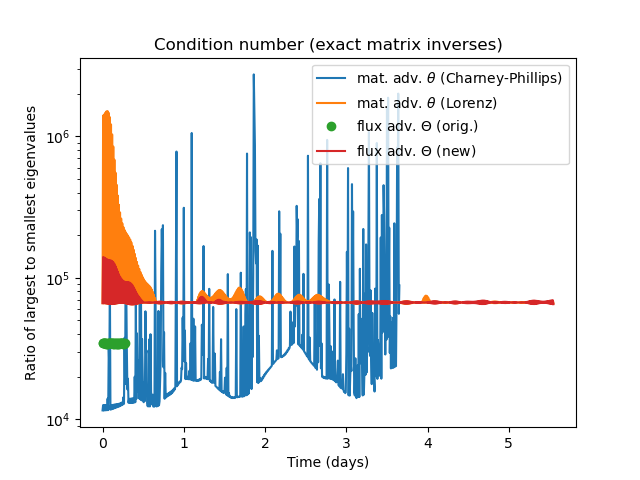}
\includegraphics[width=0.32\textwidth,height=0.24\textwidth]{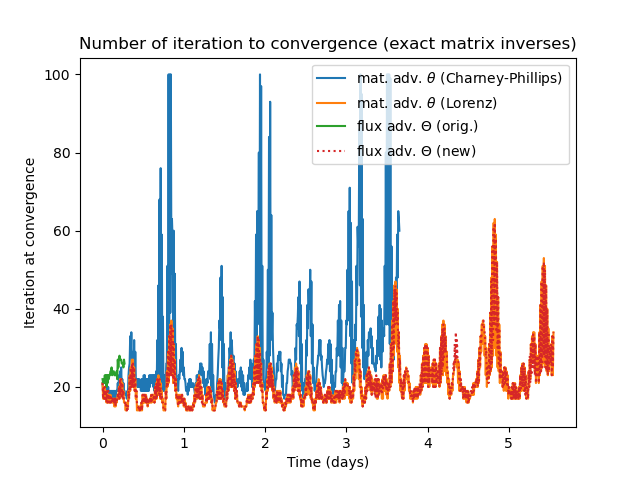}
\caption{Energy conservation error (left), Helmholtz operator condition number (center)
	and number of nonlinear iterations to convergence (right) for the three different preconditioners
	using exact matrix inverses.}
\label{fig::vert_1d_1}
\end{center}
\end{figure}

Figure~\ref{fig::vert_1d_1} shows the energy conservation error, Helmholtz matrix condition number, 
and number of nonlinear iterations to convergence for the different preconditioners using direct
matrix inverses and iterating each time step to convergence below a tolerance of $10^{-14}$ for
$|\delta\rho_h|/|\rho_h|$, $|\delta\mathcal{T}_h|/|\mathcal{T}_h|$, $|\delta\Pi_h|/|\Pi_h|$ for
$\mathcal{T}_h\in\theta_h, \Theta_h$ (we exclude the velocity residual from the convergence criteria, 
since the absolute values of this are small and somewhat volatile). While the original flux form
preconditioner detailed in Section~\ref{sec::flux_orig_prec} and \cite{Lee(2021)} fails after just 40
time steps, and the material form preconditioner fails after 527 time steps, the new flux form
and energy conserving material formulation run stably for the full duration of 
the simulation (800 steps). The instability of the Charney-Phillips material form scheme is reflected 
in the long term growth in energy and more variation in the condition number of its Helmholtz operator. 
The condition number of the Charney-Phillips grid Helmholtz operator is generally lower than that for 
the Lorenz grid operators (albeit with much more variation), suggesting faster convergence for the 
inner linear solve. However the Charney-Phillips material form also requires on average more outer nonlinear 
iterations to achieve convergence. This is perhaps a result of the smaller condition number such that 
the Charney-Phillips grid operator may not span the eigenvalues of the true 
Jacobian as effectively as the Lorenz grid operators. The new flux form and energy conserving 
material formulation are 
very similar in all cases, with perhaps the main difference being in the condition number of the
matrix at early times during the process of hydrostatic adjustment, suggesting that the new flux form
will require fewer linear solver iterations for this regime.

\begin{figure}[!hbtp]
\begin{center}
\includegraphics[width=0.32\textwidth,height=0.24\textwidth]{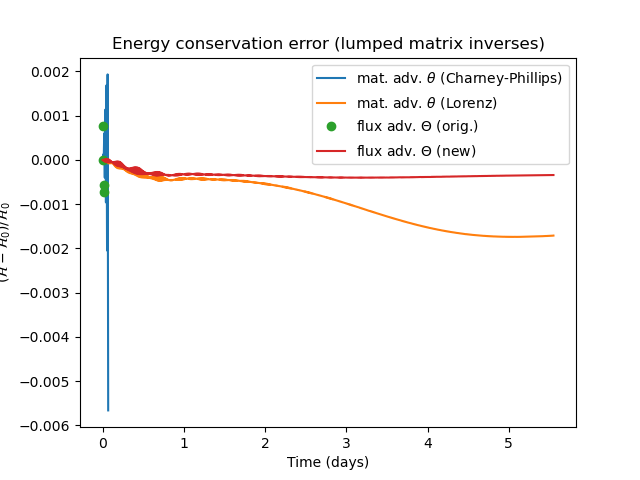}
\includegraphics[width=0.32\textwidth,height=0.24\textwidth]{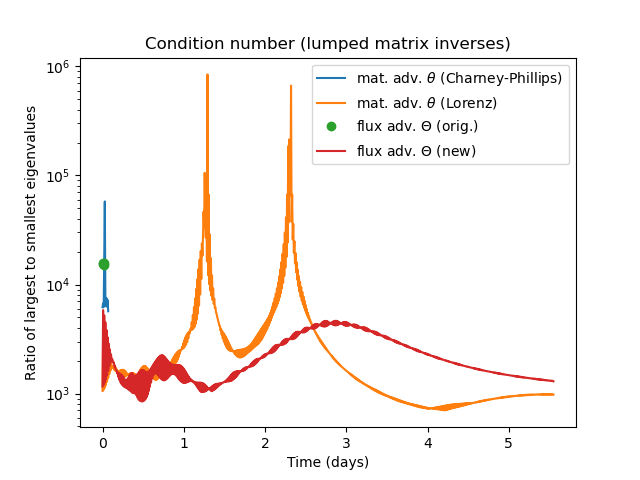}
\includegraphics[width=0.32\textwidth,height=0.24\textwidth]{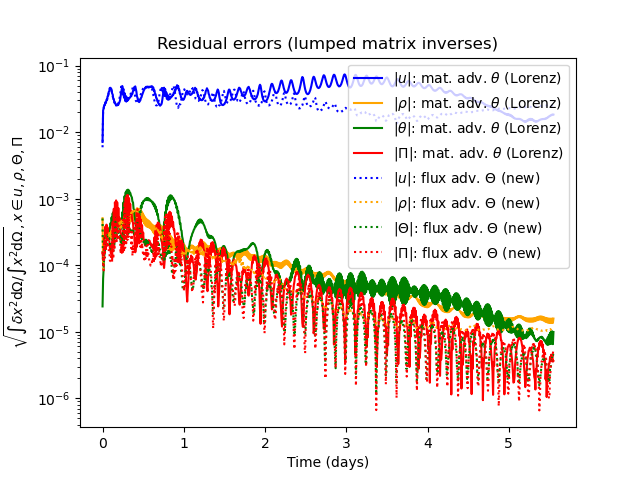}
\caption{Energy conservation error (left), Helmholtz operator condition number (center)
	and residual error at the final Newton iteration (right) 
	for the three different preconditioners
	using approximate matrix inverses and four Newton iterations per time step.}
\label{fig::vert_1d_2}
\end{center}
\end{figure}

For large scale applications on distributed memory computing architectures it is not feasible to 
directly compute matrix inverses, or to iterate to convergence, as for the configurations presented
in Figure~\ref{fig::vert_1d_1}. Consequently all intermediate matrix inverses are replaced by lumped
approximations, with the exception of the Helmholtz operator, and the resulting linear system
is solved only a finite number of times per time step rather than to nonlinear convergence. 
Results using lumped matrix inverses and four Newton nonlinear iterations are presented in Figure 
\ref{fig::vert_1d_2} for the energy conservation errors, Helmholtz matrix condition number and 
residual errors (for the two formulations that ran stably to completion) at the final Newton 
iteration respectively. For this configuration both the
material form preconditioner and the original flux form preconditioner are stable for only a small
number of time steps (11 and 4 respectively), while once again the new flux form preconditioner and
the energy conserving material formulation on the Lorenz grid are
stable for the full 800 time steps with no sign of instability in the energy conservation error, 
condition number or residual error.

Using only a finite number of nonlinear iterations the differences between the new flux form and energy conserving
material formulations are more apparent. The material formulation shows more energy damping, as well as
more variation in the matrix condition number, which is greater for moderate times, but lower for long 
times, perhaps due to the effects of the damped energy. The residual errors at the fourth Newton iteration
are similar for both formulations, but marginally lower for the new flux formulation in most cases.

Since the new flux form preconditioner is applied to the same residuals \eqref{eq::ce_disc}
as the original flux form preconditioner as presented in Section~\ref{sec::flux_orig_prec}, 
we can attribute the enhanced stability for both the converged and finite number of nonlinear iterations 
solutions to the new preconditioner. 

\subsection{2D non-hydrostatic gravity wave}

In order to verify the new preconditioner, this is applied to the
solution of a standard test case for a non-hydrostatic gravity wave travelling with an
initial mean flow of $20ms^{-1}$ in the horizontal direction \cite{SK(1994)},\cite{GR(2008)},\cite{Melvin(2019)}.
The domain is of size $[0,3\times 10^5]\times[0,10^4]m$ and is discretised using
$300\times 10$ lowest order elements and is run for a simulation time of $3000s$ 
using a time step of $\Delta t=20s$.

As for the 1D tests, no upwinding is applied such that $c=0$ in both
\eqref{eq::ce_disc_u} and \eqref{eq::ce_disc_Theta}. However we note that since this
term is applied in a skew symmetric way in these two equations, the presence of upwinding
would not change the energy conservation result presented above.

In order to suppress grid scale oscillations
we also introduce an interior penalty term, similar to the continuous interior penalty term \cite{BE(2007)} 
which penalises against jumps in the mass flux gradient at element boundaries, but with the addition of a second
term which penalises against jumps in the tangent mass flux also. This is added to the discrete momentum equation 
\eqref{eq::ce_disc_u} as
\begin{equation}\label{eq::cip}
	\int \Delta x^2u_m\ldblbrace\overline{\alpha}_h\rdblbrace
	[\![\nabla\boldsymbol{v}_h\cdot\hat{\boldsymbol{n}}_{\Gamma}]\!]
	[\![\nabla\overline{\boldsymbol{F}}_h\cdot\hat{\boldsymbol{n}}_{\Gamma}]\!] +
	u_m\ldblbrace\overline{\alpha}_h\rdblbrace
	[\![\boldsymbol{v}_h\cdot\hat{\boldsymbol{t}}_\Gamma]\!]
	[\![\overline{\boldsymbol{F}}_h\cdot\hat{\boldsymbol{t}}_\Gamma]\!]
	\mathrm{d}\Gamma
	\qquad\forall\boldsymbol{v}_h\in\mathbb{W}_2\ ,
\end{equation}
where $\Delta x$ is the element spacing, $u_m=0.5$ is damping parameter that scales with the mean velocity,
$\hat{\boldsymbol{t}}_\Gamma$ is the tangent unit normal vector along the element boundary, and 
$\overline{\alpha}_h$ is the discrete inverse density averaged over the time level, computed within each element as
\begin{equation}
	\int\phi_h\overline{\rho}_h\overline{\alpha}_h\mathrm{d}\Omega=\int\phi_h\mathrm{d}\Omega\qquad\forall\phi_h\in\mathbb{W}_3.
\end{equation}
A linearised form of this term is also added to the Helmholtz operator, via a correction to 
the operator $\boldsymbol{\mathsf{M}}_{2,R}$
as
\begin{equation}
	\Delta t\int \Delta x^2u_m
	[\![\nabla\boldsymbol{v}_h\cdot\hat{\boldsymbol{n}}_{\Gamma}]\!]
	[\![\nabla\boldsymbol{w}_h\cdot\hat{\boldsymbol{n}}_{\Gamma}]\!] + 
	u_m[\![\boldsymbol{v}_h\cdot\hat{\boldsymbol{t}}_\Gamma]\!]
	[\![\boldsymbol{w}_h\cdot\hat{\boldsymbol{t}}_\Gamma]\!]
	\mathrm{d}\Gamma
	\qquad\forall\boldsymbol{v}_h,\boldsymbol{w}_h\in\mathbb{W}_2\ .
\end{equation}

The variational derivatives of the discrete energy with respect to the prognostic variables are given by 
\eqref{eq::diag_F}, \eqref{eq::diag_Phi} and \eqref{eq::ce_disc_Pi} respectively. Setting the test functions
for the prognostic equations 
\eqref{eq::ce_disc_u}, \eqref{eq::ce_disc_rho}, \eqref{eq::ce_disc_Theta}
from these as $\overline{\boldsymbol{F}}_h$, $\overline{\Phi}_h$ and $\overline{\Pi}_h$ respectively leads
to the cancellation of all forcing terms, with the exception of the continuous interior penalty term, and
the evolution of total energy \eqref{eq::energy}
as 
\begin{equation}\label{eq::pen_term_en_diss}
	\frac{\mathrm{d}\mathcal{H}_h}{\mathrm{d}t} = 
	\int\overline{\boldsymbol{F}}_h\cdot\frac{\partial\boldsymbol{u}_h}{\partial t} + 
	\overline{\Phi}_h\frac{\partial\rho_h}{\partial t} + 
	\overline{\Pi}_h\frac{\partial\Theta_h}{\partial t}\mathrm{d}\Omega = -
	\int \Delta x^2u_m\ldblbrace\overline{\alpha}_h\rdblbrace
	([\![\nabla\overline{\boldsymbol{F}}_h\cdot\hat{\boldsymbol{n}}_{\Gamma}]\!])^2 +
	u_m\ldblbrace\overline{\alpha}_h\rdblbrace
	([\![\overline{\boldsymbol{F}}_h\cdot\hat{\boldsymbol{t}}_\Gamma]\!])^2
	\mathrm{d}\Gamma\ .
\end{equation}
Assuming $\overline{\rho}_h,\overline{\alpha}_h>0$ and $u_m\ge 0$, 
setting $\boldsymbol{v}_h=\overline{\boldsymbol{F}}_h$ in 
\eqref{eq::ce_disc_u} and \eqref{eq::cip}
results in a term that is always energy neutral 
(for smooth normal mass flux gradients and smooth tangent mass fluxes) or damping. 

In addition to the bottom and top boundary conditions described above, 
we apply periodic boundary conditions 
in the horizontal dimension. Since it is not practical to directly compute matrix inverses 
for larger problem sizes in parallel, the two dimensional configuration uses the lumped 
inverse formulation with four nonlinear iterations per time step. 

The gravity wave is driven by an initial perturbation to an otherwise hydrostatically
balanced potential temperature profile of the form
\begin{equation}
	\theta_p(t=0) = \frac{0.01\sin(\pi z/H)}{1 + (x-x_c)^2/a_c^2}\ ,
\end{equation}
where $H=10^4m$ is the domain height, $a_c=5\times10^3m$, $x_c=10^4m$ is the center of 
the perturbation, while the initial mean potential temperature is given as a function
of the Brunt-V\"ais\"al\"a frequency $\mathcal{N}=0.01s^{-1}$ as $\theta_m(t=0)=\theta_0e^{\mathcal{N}^2z/g}$,
where $\theta_0=300^{\circ}K$.
The initial Exner pressure and density profiles are given respectively as
\begin{equation}
	\Pi(t=0) = c_p + \frac{g^2(e^{-\mathcal{N}^2z/g}-1)}{\theta_0\mathcal{N}^2}\ ,
\end{equation}
and
\begin{equation}
	\rho(t=0) = \frac{p_0}{R\theta_m}\Bigg(\frac{\Pi}{c_p}\Bigg)^{c_v/R}\ ,
\end{equation}
while the initial density weighted potential temperature given as $\Theta=\rho(\theta_m+\theta_p)$.
See the appendix for a full description of the constants above.
In order to ensure an initial state of hydrostatic balance, these initial conditions 
are applied as finite element projections, rather than analytic functions.

We also run the gravity wave test case for a two dimensional configuration of the 
material transport of potential temperature on the Charney-Phillips grid as described in 
Section \ref{sec::mat_prec}. This solver uses the same spatial and temporal resolution, 
initial conditions and stabilisation term as for the new flux form transport Lorenz grid 
preconditioner as described above.

\begin{figure}[!hbtp]
\begin{center}
\includegraphics[width=0.72\textwidth,height=0.24\textwidth]{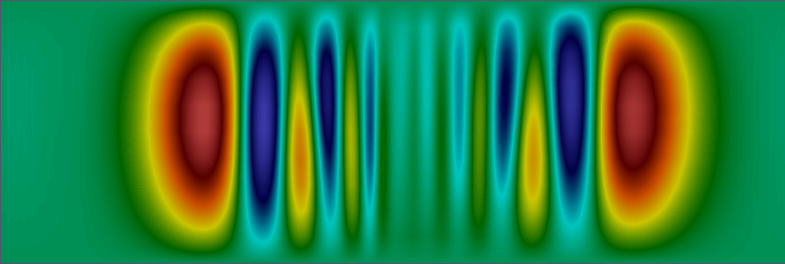}
\includegraphics[width=0.72\textwidth,height=0.24\textwidth]{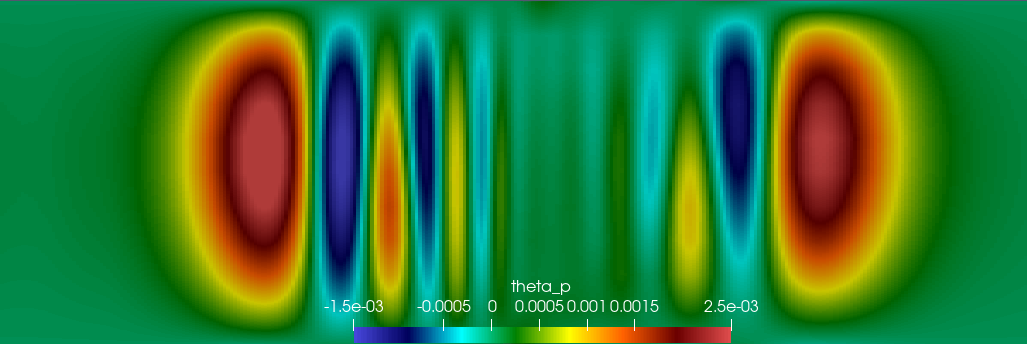}
	\caption{Potential temperature perturbation from the mean initial state, 
	$\theta_p = \theta-\theta_m(t=0)$ at time $3000s$ for the 2D non-hydrostatic gravity wave 
	test case for the new (top) and Charney-Phillips (bottom) formulations. For the new formulation 
	the potential temperature perturbation is projected from discontinuous space $\mathbb{W}_3$ to 
	continuous space $\mathbb{W}_0$. Colors range from $-0.00146^{\circ}K$ to $+0.00252^{\circ}K$. 
	Vertical axes is scaled by a factor of 10 with respect to the horizontal.}
\label{fig::gw_theta_pert_3000s}
\end{center}
\end{figure}

The perturbed potential temperature (difference between the current value and the mean
initial value, $\theta_m$) is given at time $3000s$ for the two different formulations
in Figure~\ref{fig::gw_theta_pert_3000s}. Since this is represented on the discontinuous piecewise constant space $\mathbb{W}_3$
for the flux form Lorenz formulation, for clarity this is projected into the continuous piecewise linear space
$\mathbb{W}_0$ in the top figure.
Despite the low resolution of the test configuration and long time step, this agrees well 
with a previous high resolution solution \cite{GR(2008)}, both in terms of position and 
amplitude of the disturbance.

\begin{figure}[!hbtp]
\begin{center}
\includegraphics[width=0.48\textwidth,height=0.36\textwidth]{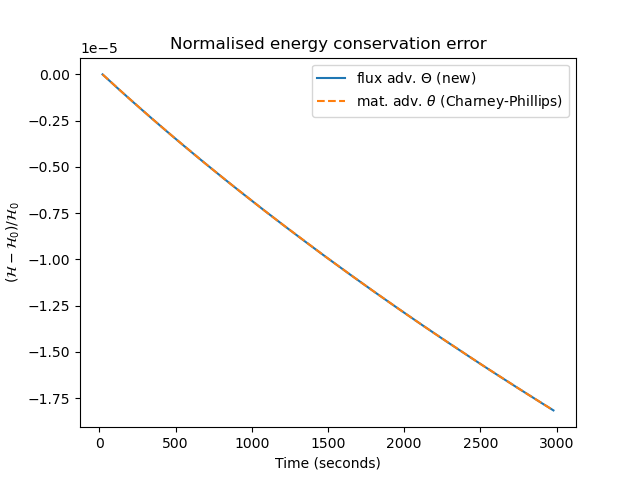}
\includegraphics[width=0.48\textwidth,height=0.36\textwidth]{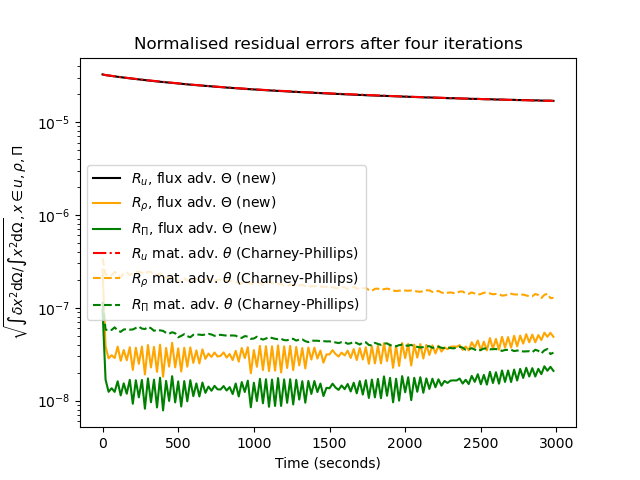}
\caption{Energy conservation error (top), and residual errors after four nonlinear 
	iterations (bottom) for the 2D non-hydrostatic gravity wave test case for the 
	new formulation and the material transport of $\theta$ on the Charney-Phillips grid.}
\label{fig::gw_cons_resid_err}
\end{center}
\end{figure}

The energy conservation and residual errors on the fourth Newton iteration as a function 
of time are given for both formulations in Figure~\ref{fig::gw_cons_resid_err}. 
There is a steady decay of energy, consistent with the application of the penalty term as 
described in \eqref{eq::pen_term_en_diss}, which is almost identical for the two formulations.
The residual errors for $\mathcal{R}_u$, $\mathcal{R}_{\rho}$ and $\mathcal{R}_{\Pi}$ are smaller 
for the new preconditioner, but increase somewhat as the solution evolves, while those for the
Charney-Phillips formulation decrease. For a larger time step of $\Delta t=30s$ at the same spatial 
resolution grid scale instabilities begin to develop for both formulations (not shown), 
but these are more pronounced for the new flux form transport Lorenz formulation.

The residuals exhibit an oscillation on a time scale of $2\Delta t$, perhaps due to the presence of acoustic modes, 
which are not explicitly resolved for a spatial resolution of $\Delta x=1000m$ and a time 
step of $\Delta t=20s$. Again, this oscillation is more pronounced for the new flux form transport 
Lorenz grid formulation. These oscillations are reduced for both formulations for a doubling of 
the spatial and temporal resolutions (not shown). 
The more pronounced nature of these oscillations for the new flux form Lorenz grid formulation is
perhaps due to the presense of the computational mode on the Lorenz grid.

The material form Charney-Phillips formulation studied here has the same spatial
discretisation and approximate Jacobian as that used in the {\tt LFRic} model \cite{Melvin(2019)},\cite{Maynard(2020)}.
However the {\tt LFRic} model uses upwinded finite volume transport for the potential temperature, 
and optionally for the velocity as well. This upwinding provides additional stabilisation 
such that the stabilisation term presented here \eqref{eq::cip} is not required in that model.

Execution times are identical for both formulations; 655 seconds on 6 processors with a 
full node reserved for the computation. However the new flux form Lorenz grid preconditioner
requires slightly fewer iterations for each linear GMRES solve at each Newton iteration, averaging
50.58 iterations per linear solve for 52.44 iterations per linear solve for the material transport
of potential temperature on the Charney-Phillips grid.

\subsection{Density current}

The new preconditioner is further verified against a standard test case for a sinking cold 
bubble that falls to the bottom boundary before progressing outwards as a density current
\cite{Straka(1993)}, \cite{GR(2008)}, \cite{Melvin(2019)}. The 2D domain is configured with 
periodic boundary conditions in the horizontal as $[-2.56\times 10^4,+2.56\times 10^4]\times[0,6.4\times 10^3]m$ using 
$864\times108$ regularly spaced lowest order elements for a resolution of $\Delta x\approx 59.26m$ and a
time step of $\Delta t=2.5s$. The initial conditions are similar to those for the gravity 
wave test case above, with an initial state of isothermal ($\theta_0=300^{\circ}K$) 
hydrostatic balance overlaid with a potential temperature perturbation, in this 
case a cold bubble of the form
\begin{equation}
	\theta_p(t=0) = \theta_0 - 7.5(1 + \cos(\pi r))\ ,
\end{equation}
where $r = \sqrt{(x/4000)^2 + ((z-3000)/2000)^2}$ for $r < 1$, and the initial Exner pressure
and density given respectively as 
\begin{equation}
	\Pi(t=0) = c_p\Big(1-\frac{gz}{c_p\theta_0}\Big)\ ,
\end{equation}
\begin{equation}
	\rho(t=0) = \frac{p_0}{R\theta_0}\Bigg(\frac{\Pi}{c_p}\Bigg)^{c_v/R}.
\end{equation}

Unlike the previous test case, the density current test case calls for a specific viscous term with a
coefficient of $\nu=75.0m^2/s$ to be applied to both the momentum \eqref{eq::ce_disc_u} and potential
temperature transport \eqref{eq::ce_disc_Theta} equations. These are implemented (for lowest order elements) respectively as
\begin{equation}
	\int\nu\nabla\boldsymbol{v}_h:\nabla\boldsymbol{\overline{u}}_h\mathrm{d}\Omega + 
	\int
	\frac{\nu}{\Delta x}
	[\![\boldsymbol{v}_h\cdot\hat{\boldsymbol{t}}_{\Gamma}]\!]
	[\![\overline{\boldsymbol{u}}_h\cdot\hat{\boldsymbol{t}}_{\Gamma}]\!] +
	\nu\Delta x
	[\![\nabla\boldsymbol{v}_h\cdot\hat{\boldsymbol{n}}_{\Gamma}]\!]\cdot
	[\![\nabla\overline{\boldsymbol{u}}_h\cdot\hat{\boldsymbol{n}}_{\Gamma}]\!]
	\mathrm{d}\Gamma
	\qquad\forall\boldsymbol{v}_h\in\mathbb{W}_2\ ,
\end{equation}
\begin{equation}
	\int\frac{\nu}{\Delta x}\ldblbrace\overline{\rho}_h\rdblbrace
	[\![\psi_h]\!]
	[\![\overline{\theta}_h]\!]\mathrm{d}\Gamma\qquad\forall\psi_h\in\mathbb{W}_3.
\end{equation}
These viscous terms will introduce small diffusive timescales to the problem. To ensure
these do not have to be explicitly resolved by the time stepping scheme, linearised versions
of these terms are also added to the respectively diagonal blocks of the approximate Jacobian
\eqref{eq::J} as
\begin{equation}
	\Delta t\int \nu\nabla\boldsymbol{v}_h:\nabla\boldsymbol{w}_h\mathrm{d}\Omega +
	\Delta t\int \frac{\nu}{\Delta x}
	[\![\boldsymbol{v}_h\cdot\hat{\boldsymbol{t}}_{\Gamma}]\!]
	[\![\boldsymbol{w}_h\cdot\hat{\boldsymbol{t}}_{\Gamma}]\!] +
	\nu\Delta x
	[\![\nabla\boldsymbol{v}_h\cdot\hat{\boldsymbol{n}}_{\Gamma}]\!]\cdot
	[\![\nabla\boldsymbol{w}_h\cdot\hat{\boldsymbol{n}}_{\Gamma}]\!]
	\mathrm{d}\Gamma
	\qquad\forall\boldsymbol{v}_h,\boldsymbol{w}_h\in\mathbb{W}_2\ ,
\end{equation}
\begin{equation}\label{eq::visc_eta}
	\Delta t\int\frac{\nu}{\Delta x}
	[\![\psi_h]\!]
	[\![\phi_h]\!]\mathrm{d}\Gamma\qquad\forall\psi_h,\phi_h\in\mathbb{W}_3.
\end{equation}
In the case of the linearised form of the potential temperature operator \eqref{eq::visc_eta},
this omits several additional terms that arise from the transformation to an entropy
residual equation \eqref{eq::R_eta}, however it seems to account for sufficient stiffness
to allow for robust convergence. 

Note that we have not made comparisons to the material form Charney-Phillips grid formulation
for this test case, since the viscous term for the potential temperature described above would need
to be implemented differently, owing to the difference in finite element spaces used to represent 
the Charney-Phillips space, such that the comparison would be somewhat indirect.

As before, no upwinding of the potential temperature is applied ($c=0$).
While upwinding or some other form of advective stabilisation is often applied to the potential 
temperature transport equation \cite{GR(2008)},\cite{Melvin(2019)},\cite{Bendall(2020)},\cite{LP(2021)}, 
this would lead to excessively diffusive solutions for our piecewise constant representation of the 
density weighted potential temperature, since this would result
in potential temperature fluxes that were only first order accurate. One possible means of negating this issue
could be to project the potential temperature into a higher order space before applying upwinded flux,
as has been done previously for a low order mixed finite element formulation \cite{Bendall(2019)}.

\begin{figure}[!hbtp]
\begin{center}
\includegraphics[width=0.72\textwidth,height=0.18\textwidth]{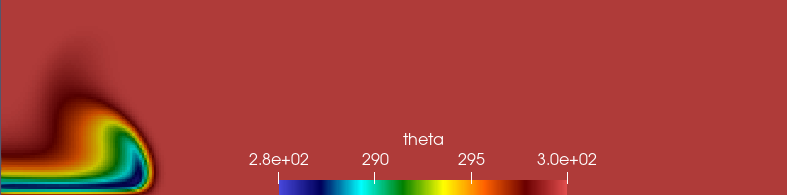}
\includegraphics[width=0.72\textwidth,height=0.18\textwidth]{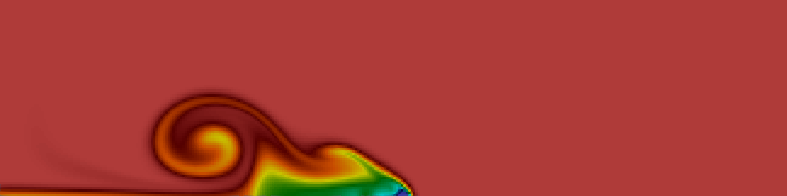}
\includegraphics[width=0.72\textwidth,height=0.18\textwidth]{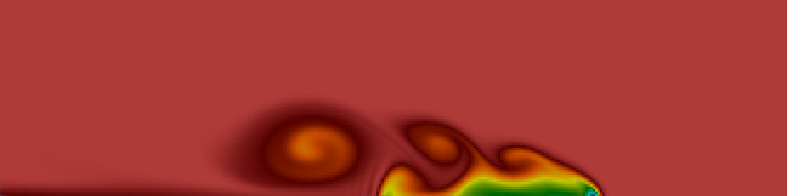}
	\caption{Potential temperature for the density current test case at 
	$300s$, $600s$, $900s$ over the subdomain of $[0,1.92\times 10^4]\times[0,4.8\times 10^3]m$.
	Colors range from $285^{\circ}K$ to $300^{\circ}K$.}
\label{fig::theta_straka}
\end{center}
\end{figure}

\begin{figure}[!hbtp]
\begin{center}
\includegraphics[width=0.48\textwidth,height=0.36\textwidth]{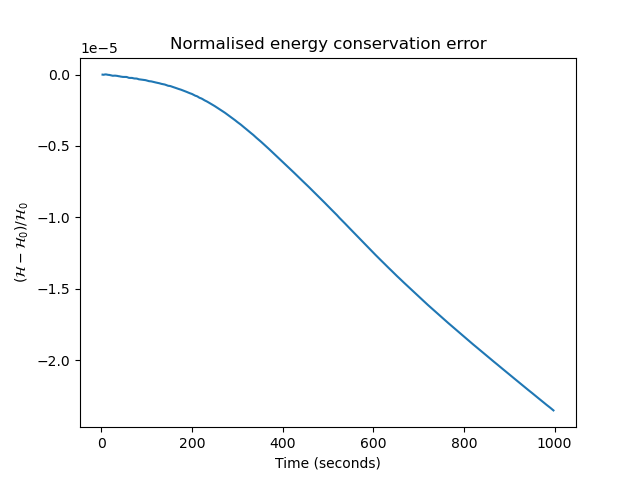}
\includegraphics[width=0.48\textwidth,height=0.36\textwidth]{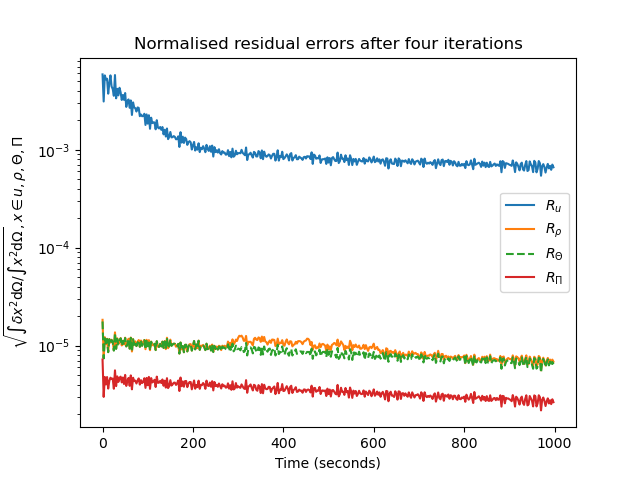}
\caption{Energy conservation error (top), and residual errors after four nonlinear 
	iterations (bottom) for the 2D density current test case.}
\label{fig::straka_cons_resid_err}
\end{center}
\end{figure}

The potential temperature, $\theta_h$ for the sinking bubble is given for times $300s$, $600s$
and $900s$ for the sub-domain $[0,1.92\times 10^4]\times[0,4.8\times 10^3]m$ in Fig. \ref{fig::theta_straka},
with colors ranging from $285^{\circ}K$ and $300^{\circ}K$. These results match well against
previously published results
\cite{Straka(1993)}, \cite{GR(2008)}, \cite{Melvin(2019)}, in terms of position and shape,
with three rotors clearly visible at the final time. The energy conservation and residual errors
at the fourth (final) nonlinear iteration are given as a function of time in Fig. \ref{fig::straka_cons_resid_err}.
While the viscous term given here is not provably energy dissipating, as is the case for the penalty
term applied for the gravity wave test case, the energy decays monotonically with a similar magnitude.
The residual errors at the final Newton iteration reduce with timestep. As with the gravity wave test, these residual errors
also show a small temporal oscillation, potentially due to unresolved acoustic modes. One potential means of
suppressing this oscillation may be the use of an L-stable time integration scheme, such as a Rosenbrock-Wanner
method, which has been previously shown to suppress small scale temporal oscillations that are present for
time centered integration schemes for compressible atmospheric simulations using Helmholtz preconditioning 
\cite{Lee(2024)}.

\section{Conclusions}
\label{sec::conclusions}

A new Helmholtz preconditioner for the compressible Euler equations with flux form density weighted
potential temperature transport on the Lorenz grid using mixed finite elements is presented. A 
transformation of the density and density weighted potential temperature residuals into a material
transport equation for the thermodynamic entropy allows the new preconditioner to preserve a 
similar block structure as an existing preconditioner for the material transport of potential 
temperature on the Charney-Phillips grid. One dimensional comparisons against this and another 
previous formulation for material form transport on the Lorenz grid show the new preconditioner
to be more stable, both when run to convergence using exact matrix inverses and for when using
just four Newton nonlinear iterations per time step with approximate lumped matrix inverses within the
compound Helmholtz operator. Comparisons to an alternative energy conserving formulation for
material form potential temperature transport on the Lorenz grid show little variation for 
converged solutions, but less energy damping and lower condition numbers for a finite number of
nonlinear iterations.

The new preconditioner is further verified for existing two dimensional test cases within vertical
slice configurations. Comparisons to an existing formulation using material form transport of
potential temperature on the Charney-Phillips grid in two dimensions shows smaller residual errors 
and fewer iterations to converge for each inner linear GMRES solve, but a higher propensity to unstable
oscillations with increased time step.

In the present case the preconditioner is applied using a time centered implicit 
discretisation, which can give rise to small scale temporal oscillations due to unresolved 
modes, as observed in the residual errors for the test cases presented here. One means of suppressing 
these oscillations and increasing model stability may be the application the preconditioner to an L-stable Rosenbrock-Wanner
time integration scheme, which as been previously observed to suppress such modes for compressible 
atmospheric simulations \cite{Lee(2024)}.

Another potential future avenue of research may be transform the residual for the material transport equation for $\theta$
on the Charney-Phillips grid into a residual expression for $\eta$. This would allow for a Jacobian linearisation similar 
to that for the new flux form Lorenz grid preconditioner, but on the Charney-Phillips grid. This might potentially lead to
improved conditioning and convergence for that that existing formulation.

\section*{Acknowledgements}
This work was (partially) supported by computational resources provided by the Australian Government through the
National Computational Infrastructure (NCI) under the ANU Merit Allocation Scheme. Kieran Ricardo would like to acknowledge 
the Australian Government through the Australian Government Research Training
Program (RTP) Scholarship, and the Bureau of Meteorology through research contract KR2326. 
We would also like to thank Dr. Matthew Wheeler for his helpful comments on an early draft of this manuscript.

\section*{Data availability}
Source code and model output data will be made available upon request.

\section*{Appendix: Vertical profiles for the hydrostatically balanced 1D atmosphere}

The 1D vertical profiles are taken from \cite{UMJS(2014)}. Constants are given as:
$T_e=310K$, $T_p=240K$, $g=9.80616ms^{-2}$, $\gamma=0.005Km^{-1}$, $p_0=10^{5}Pa$, $c_p=1004.5Jkg^{-1}K^{-1}$, $R=287Jkg^{-1}K^{-1}$.
In order to determine the reference profiles, we first introduce some intermediate variables as:
\begin{subequations}
	\begin{align}
		T_0 &= (T_e+T_p)/2 \\
		B &= (T_e-T_p)/((T_e+T_p)T_p) \\
		C &= 5(T_e-T_p)/(2T_eT_p) \\
		D &= \cos^3(2\pi/9) - (3/5)\cos^5(2\pi/9) \\
		E &= (gz/(2RT_0))^2 \\
		\tau_1 &= T_0^{-1}e^{\gamma z/T_0} + B(1-2E)e^{-E}\\
		\tau_2 &= C(1-2E)e^{-E} \\
		\chi_1 &= \gamma^{-1}(e^{\gamma z/T_0}-1) + Bze^{-E} \\
		\chi_2 &= Cze^{-E}\ .
	\end{align}
\end{subequations}
From the intermediate values given above, the initial vertical profiles for the
temperature, pressure, Exner pressure, density and potential temperature may be
given respectively as
\begin{subequations}
	\begin{align}
		T &= (\tau_1 - \tau_2D)^{-1} \\
		p &= p_0e^{-g\chi_1/R + g\chi_2D/R} \\
		\Pi &= c_p(p/p_0)^{R/c_p} \\
		\rho &= p/(RT) \\
		\theta &= T(p_0/p)^{R/c_p}\ .
	\end{align}
\end{subequations}


\begin{thebibliography}{32}
\providecommand{\natexlab}[1]{#1}
\providecommand{\url}[1]{\texttt{#1}}
\expandafter\ifx\csname urlstyle\endcsname\relax
  \providecommand{\doi}[1]{doi: #1}\else
  \providecommand{\doi}{doi: \begingroup \urlstyle{rm}\Url}\fi

\bibitem[(Arakawa and Connor, 1996)]{Arakawa(1996)}
A~Arakawa and C~S~Connor.
\newblock Vertical differencing of the primitive equations based on the 
  {C}harney-{P}hillips grid in hybrid $\sigma-p$ vertical coordinates.
\newblock \emph{Mon. Wea. Rev.}, 124:\penalty0 511--528, 1996.

\bibitem[(Badia and Verdugo, 2020)]{Badia(2020)}
S.~Badia and F.~Verdugo. 
\newblock {G}ridap: An extensible finite element toolbox in {J}ulia. 
\newblock \emph{J. Open Source Softw.}, 5:\penalty0 2520, 2020.

\bibitem[(Badia et~al. 2022)]{Badia(2022)}
S.~Badia, A.~F.~Mart\'{\i}n, and F.~Verdugo.
\newblock {G}ridap{D}istributed: a massively parallel finite element toolbox in {J}ulia.
\newblock \emph{J. Open Source Softw.}, 7:\penalty0 4157, 2022.

\bibitem[(Bauer and Cotter, 2018)]{Bauer(2018)}
W~Bauer and C~J~Cotter.
\newblock Energy–enstrophy conserving compatible finite element schemes 
  for the rotating shallow water equations with slip boundary conditions.
\newblock \emph{J. Comp. Phys.}, 373:\penalty0 171--187, 2018.

\bibitem[(Bendall et~al. 2019)]{Bendall(2019)}
T~M~Bendall, C~J~Cotter, J~Shipton.
\newblock The 'recovered space' advection scheme for lowest-order compatible
  finite element methods.
\newblock \emph{J. Comp. Phys.}, 390:\penalty0 342--358, 2019.

\bibitem[(Bendall et~al. 2020)]{Bendall(2020)}
T~M~Bendall, T~H~Gibson, J~Shipton, C~J~Cotter, B~Shipway.
\newblock A compatible finite-element discretisation for the moist compressible 
  {E}uler equations.
\newblock \emph{Q. J. R. Meteorol. Soc.}, 146:\penalty0 3187--3205, 2020.

\bibitem[(Bendall et~al. 2023)]{Bendall(2023)}
T~M~Bendall, N~Wood, J~Thuburn, C~J~Cotter.
\newblock A solution to the trilemma of the moist {C}harney–{P}hillips staggering.
\newblock \emph{Q. J. R. Meteorol. Soc.}, 149:\penalty0 262--276, 2023.

\bibitem[(Betteridge et~al. 2023)]{Betteridge(2023)}
J~D~Betteridge, C~J~Cotter, T~H~Gibson, M~J~Griffiths, T~Melvin and E~H~M\"uller.
\newblock Hybridised multigrid preconditioners for a compatible finite-element 
  dynamical core.
\newblock \emph{Q. J. R. Meteorol. Soc.}, 149:\penalty0 2454--2476, 2023.

\bibitem[(Burman and Ern, 2007)]{BE(2007)}
E~Burman and A~Ern.
\newblock Continuous interior penalty $hp$-finite element methods for advection
  and advection-diffusion equations.
\newblock \emph{Mathematics of Computation}, 76:\penalty0 1119--1140, 2007.

\bibitem[(Cohen and Hairer, 2011)]{CH(2011)}
D~Cohen and E~Hairer.
\newblock Linear energy-preserving integrators for {P}oisson systems.
\newblock \emph{BIT Numer. Math.}, 51:\penalty0 91--101, 2011.

\bibitem[(Eldred et~al. 2019)]{Eldred(2019)}
C~Eldred, T~Dubos, and E~Kritsikis.
\newblock A quasi-{H}amiltonian discretization of the thermal shallow
  water equations.
\newblock \emph{J. Comp. Phys.}, 379:\penalty0 1--31, 2019.

\bibitem[(Gibson et~al, 2020)]{Gibson(2020)}
T~H~Gibson, L~Mitchell, D~A~Ham, C~J~Cotter.
\newblock Slate: extending {F}iredrake's domain-specific abstraction to
  hybridized solvers for geoscience and beyond.
\newblock \emph{Geosci. Model. Dev.}, 13:\penalty0 735--761, 2020.

\bibitem[(Giraldo and Restelli, 2008)]{GR(2008)}
F~X~Giraldo and M~Restelli.
\newblock A study of spectral element and discontinuous {G}alerkin methods 
  for the Navier–Stokes equations in nonhydrostatic mesoscale atmospheric 
  modeling: {E}quation sets and test cases.
\newblock \emph{J. Comp. Phys.}, 227:\penalty0 3849--3877, 2008.

\bibitem[(Lee, 2021)]{Lee(2021)}
D~Lee.
\newblock An energetically balanced, quasi-{N}ewton integrator for
  non-hydrostatic vertical atmospheric dynamics.
\newblock \emph{J. Comp. Phys.}, 429:\penalty0 109988, 2021.

\bibitem[(Lee and Palha, 2021)]{LP(2021)}
D~Lee and A~Palha.
\newblock Exact spatial and temporal balance of energy exchanges within a 
  horizontally explicit/vertically implicit non-hydrostatic atmosphere.
\newblock \emph{J. Comp. Phys.}, 440:\penalty0 110432, 2021.

\bibitem[(Lee, 2024)]{Lee(2024)}
D~Lee.
\newblock A comparison of Rosenbrock–Wanner and Crank–Nicolson time 
  integrators for atmospheric modelling
\newblock \emph{Q. J. Royal Meteorol. Soc.}, 758:\penalty0 462--483, 2024.

\bibitem[(Maynard et~al. 2020)Maynard, Melvin, and M\"uller]{Maynard(2020)}
C~Maynard, T~Melvin, and E~H M\"uller.
\newblock Multigrid preconditioners for the mixed finite element dynamical core
  of the {LFR}ic atmospheric model.
\newblock \emph{Q. J. R. Meteorol. Soc.}, 146:\penalty0 3917--3936, 2020.

\bibitem[(Melvin et~al. 2018)Melvin, Benacchio, Thuburn, and Cotter]{Melvin(2018)}
T~Melvin, T~Benacchio, J~Thuburn, and C~Cotter.
\newblock Choice of function spaces for thermodynamic variables in mixed
  finite-element methods
\newblock \emph{Q. J. R. Meteorol. Soc.}, 144:\penalty0 900--916, 2018.

\bibitem[(Melvin et~al. 2019)Melvin, Benacchio, Shipway, Wood, Thuburn, and
  Cotter]{Melvin(2019)}
T~Melvin, T~Benacchio, B~Shipway, N~Wood, J~Thuburn, and C~Cotter.
\newblock A mixed finite-element, finite-volume, semi-implicit discretisation
  for atmospheric dynamics: {C}artesian geometry.
\newblock \emph{Q. J. R. Meteorol. Soc.}, 145:\penalty0 1--19, 2019.

\bibitem[(Natale et~al. 2016)]{Natale(2016)}
A~Natale, J~Shipton, and C~J~Cotter.
\newblock Compatible ﬁnite element spaces for geophysical ﬂuid dynamics.
\newblock \emph{Dyn. Stat. Climate Sys.}, 1:\penalty0 1--31, 2016.

\bibitem[(Reddy et~al. 2023)]{Reddy(2023)}
S~Reddy, M~Waruszewski, F~A~V~de Braganca Alves, and F~X~Giraldo.
\newblock Schur complement {IM}plicit-{EX}plicit formulations for 
  discontinuous {G}alerkin non-hydrostatic atmospheric models.
\newblock \emph{J. Comp. Phys.}, 491:\penalty0 112361, 2023.

\bibitem[(Ricardo et~al. 2024)Ricardo, Lee, Duru]{Ricardo(2024)}
K~Ricardo, D~Lee, and K~Duru.
\newblock Entropy and energy conservation for thermal atmospheric dynamics 
  using mixed compatible finite elements
\newblock \emph{J. Comp. Phys.}, 496:\penalty0 112605, 2024.

\bibitem[(Skamarock and Klemp, 1994)]{SK(1994)}
W~C~Skamarock and J~B~Klemp.
\newblock Efficiency and Accuracy of the {K}lemp-{W}ilhelmson Time-Splitting
  Technique.
\newblock \emph{Mon. Wea. Rev.}, 122:\penalty0 2623--2630, 1994.

\bibitem[(Straka et~al. 1993)]{Straka(1993)}
J~M~Straka, R~B~Wilhelmson, L~J~Wicker, J~R~Anderson, K~K~Droegemeier.
\newblock Numerical solutions of a non-linear density current: a benchmark 
  solution and comparisons
\newblock \emph{Int. J. Numer. Meth. Fluids}, 17:\penalty0, 1--22, 1993.

\bibitem[(Taylor et~al. 2020)]{Taylor(2020)}
M~A~Taylor, O~Guba, A~Steyer, P~A~Ullrich, D~M~Hall, and C~Eldred.
\newblock An Energy Consistent Discretization of the Nonhydrostatic 
  Equations in Primitive Variables.
\newblock \emph{Journal of Advances in Modelling Earth Systems}, 
  12:\penalty0 (1), 2020.

\bibitem[(Thuburn and Woolings, 2005)]{TW(2005)}
J~Thuburn and T~J~Woolings.
\newblock Vertical discretizations for compressible {E}uler equation 
  atmospheric models giving optimal representation of normal modes
\newblock \emph{J. Comp. Phys.}, 203:\penalty0 386--404, 2005.

\bibitem[(Thuburn, 2022)]{Thuburn(2022)}
J~Thuburn.
\newblock Numerical entropy conservation without sacrificing 
  {C}harney–{P}hillips grid optimal wave propagation.
\newblock \emph{Q. J. R. Meteorol. Soc.}, 148:\penalty0 2755--2768, 2022.

\bibitem[(Ullrich et~al. 2014)Ullrich, Melvin, Jablonowski, and
  Staniforth]{UMJS(2014)}
P~A Ullrich, T~Melvin, C~Jablonowski, and A~Staniforth.
\newblock A proposed baroclinic wave test case for deep‐ and
  shallow‐atmosphere dynamical cores.
\newblock \emph{Q. J. R. Meteorol. Soc.}, 140:\penalty0 1590--1602, 2014.

\bibitem[(Wood et~al. 2014)Wood, Staniforth, White, Allen, Diamantakis, Gross, Melvin, Smith, Vosper, Zerroukat, Thuburn]{Wood(2014)}
N~Wood, A~Staniforth, A~White, T~Allen, M~Diamantakis, M~Gross, T~Melvin, C~Smith, S~Vosper, M~Zerroukat and J~Thuburn.
\newblock An inherently mass-conserving semi-implicit semi-{L}agrangian discretization of the deep-atmosphere global 
  non-hydrostatic equations.
\newblock \emph{Q. J. R. Meteorol. Soc.}, 140:\penalty0 1505--1520, 2014.

\bibitem[(Yeh et~al. 2002)J~C\^ot\'e, S~Gravel, A~M\'ethot, A~Patoine, M~Roch, A~Staniforth]{Yeh(2002)}
K-S~Yeh, J~C\^ot\'e, S~Gravel, A~M\'ethot, A~Patoine, M~Roch and A~Staniforth.
\newblock The {CMC}-{MRB} {G}lobal {E}nvironmental {M}ultiscale ({GEM}) Model.
  {P}art {III}: Nonhydrostatic Formulation.
\newblock \emph{Mon. Wea. Rev.}, 130:\penalty0 339--356, 2002.

\end{thebibliography}
\end{document}